\documentclass[a4paper,12pt]{amsart}

\usepackage[utf8]{inputenc}
\usepackage[T1]{fontenc}
\usepackage{lmodern}

\usepackage[margin=28mm,headheight=14pt]{geometry}
\usepackage{amsmath,amssymb,mathrsfs,mathtools}
\usepackage{tikz-cd}
\usepackage{microtype,needspace}
\usepackage{enumitem}
\usepackage[hidelinks]{hyperref}

\hypersetup{pdftitle={Equivariant extension of I-N-finite tempered distributions on real flag varieties},
 pdfsubject={Equivariant extension on real flag varieties}, bookmarksdepth=2}

\numberwithin{equation}{section}
\allowdisplaybreaks[2]
\newtheorem{theorem}{Theorem}[section]
\newtheorem{proposition}[theorem]{Proposition}
\newtheorem{lemma}[theorem]{Lemma}

\newtheorem{maintheorem}{Theorem}

\theoremstyle{definition}
\newtheorem{definition}[theorem]{Definition}
\theoremstyle{remark}
\newtheorem{remark}[theorem]{Remark}

\newcommand{\G}{\mathbf G}
\newcommand{\Q}{\mathbf Q}
\newcommand{\pP}{\mathbf P}
\newcommand{\lL}{\mathbf L}
\newcommand{\U}{\mathbf U}
\newcommand{\X}{\mathbf X}
\newcommand{\Y}{\mathbf Y}
\newcommand{\F}{\mathbf F}
\newcommand{\D}{\mathbf D}
\newcommand{\B}{\mathbf B}
\newcommand{\V}{\mathbf V}
\newcommand{\R}{\mathbb R}
\newcommand{\C}{\mathbb C}
\newcommand{\fg}{\mathfrak g}
\newcommand{\fq}{\mathfrak q}
\newcommand{\fp}{\mathfrak p}
\newcommand{\fl}{\mathfrak l}
\newcommand{\fn}{\mathfrak n}
\newcommand{\fu}{\mathfrak u}
\newcommand{\fa}{\mathfrak a}
\newcommand{\fc}{\mathfrak c}
\newcommand{\IN}{I_N}

\newcommand{\mcU}{\mathcal U}
\newcommand{\mcS}{\mathcal S}
\newcommand{\mcE}{\mathcal E}
\newcommand{\mcP}{\mathcal P}
\newcommand{\mcH}{\mathcal H}
\newcommand{\E}{\mathsf E}
\newcommand{\Ad}{\operatorname{Ad}}
\newcommand{\Lie}{\operatorname{Lie}}
\newcommand{\Sym}{\operatorname{Sym}}
\newcommand{\supp}{\operatorname{supp}}
\newcommand{\ord}{\operatorname{ord}}
\newcommand{\pr}{\operatorname{pr}}
\newcommand{\dd}{\mathrm d}
\newcommand{\id}{\mathrm{id}}
\newcommand{\Span}{\operatorname{Span}}

\newcommand{\tTheta}{\boldsymbol\Theta}
\newcommand{\sS}{\mathbf S}

\newcommand{\zZ}{\mathbf Z}
\newcommand{\rR}{\mathbf R}
\newcommand{\oO}{\mathbf O}
\newcommand{\oOm}{\boldsymbol\Omega}
\newcommand{\cC}{\mathbf C}
\newcommand{\lLd}{\boldsymbol\Lambda}

\newcommand{\ad}{\operatorname{ad}}
\newcommand{\Ht}{\operatorname{ht}}

\newcommand{\mcL}{\mathcal L}

\title[Equivariant extension on real flag varieties]
{Equivariant extension of $I_N$-finite tempered distributions on real flag varieties}

\author[J.Y.~Xiao]{Junyan Xiao}
\address{School of Mathematical Sciences, Peking University, No. 5 Yiheyuan Road, Beijing 100871, P.R. China}
\email{jyxiao25@stu.pku.edu.cn}
\date{}


\begin{document}

\begin{abstract}
Let $G$ be a real linear reductive group and let $P=L_PU_P$ and $Q=LN$ be parabolic subgroups. For a Casselman--Wallach representation $F$ of $L_P$, let $\mcE_F$ be the unnormalized homogeneous tempered bundle on $G/P$ induced from $F$. We prove that every $I_N$-finite tempered distribution on a $Q$-stable Nash open subset of $G/P$ extends to an $I_N$-finite tempered distribution on $G/P$.
\end{abstract}

\maketitle
\tableofcontents

\section{Introduction}\label{sec:intro}

\subsection{Main results}\label{ssec:result}

Let $G$ be a real reductive group as in Section~\ref{sec:conventions}, and let $P=L_PU_P$ and $Q=LN$ be parabolic subgroups, with $\fn=\Lie N$. For a Casselman--Wallach representation $F$ of $L_P$, inflated trivially across $U_P$, let
\[
 X=G/P,\qquad \mcE_F=G\times^P F
\]
be the unnormalized homogeneous tempered Fr\'echet bundle. Let $\IN\subset\mcU(\fn_\C)$ be the augmentation ideal. For a $Q$-stable Nash open subset $W\subset X$, write $\mcS'(W,\mcE_F)=\mcS(W,\mcE_F)'_b$ for the space of tempered distributions and set
\[
 \mcS'(W,\mcE_F)^{[\fn]}
 :=\{T\in\mcS'(W,\mcE_F):\IN^kT=0\text{ for some }k\ge1\}.
\]
These are the $\IN$-finite tempered distributions, or generalized $\fn$-invariants.

\Needspace{10\baselineskip}
\begin{maintheorem}\label{thm:general-extension}
Let $W\subset X$ be a $Q$-stable Nash open subset.
\begin{enumerate}[label=\textup{(\roman*)}]
\item For every Casselman--Wallach representation $F$ of $L_P$, restriction
\[
 \mcS'(X,\mcE_F)^{[\fn]}
 \longrightarrow\mcS'(W,\mcE_F|_W)^{[\fn]}
\]
is surjective.
\item The functor
\[
 F\longmapsto\mcS'(W,\mcE_F|_W)^{[\fn]}
\]
from $\mathrm{CW}(L_P)$ to complex vector spaces is contravariantly exact.
\end{enumerate}
\end{maintheorem}

We also prove lifting with support (Proposition~\ref{prop:coefficient-supported}). If $F_1\hookrightarrow F_2$ is an embedding in $\mathrm{CW}(L_P)$ and $O$ is a $Q$-orbit closed in a Nash open subset $U\subset X$, then
\[
 \mcS'_O(U,\mcE_{F_2})^{[\fn]}\longrightarrow
 \mcS'_O(U,\mcE_{F_1})^{[\fn]}
\]
is surjective. Here subscripts denote support conditions.

This lifting result and Casselman's subrepresentation theorem reduce Theorem~\ref{thm:general-extension}(i) to $P$ minimal. Fix a minimal parabolic $P_0=M_0A_0N_0$ and a standard parabolic $Q=MAN\supset P_0$. For a finite-dimensional representation $\sigma$ of $M_0$ and $\lambda\in\fa_{0,\C}^*$, put
\[
 X_0=G/P_0,\qquad \mcE_\lambda=G\times^{P_0}V_{\sigma,\lambda},
\]
where $V_{\sigma,\lambda}=\sigma\otimes e^\lambda$ and $N_0$ acts trivially. We first prove the following special case:

\begin{maintheorem}\label{thm:minimal-extension}
Let $V\subset X_0$ be a $Q$-stable Nash open subset. For every $\lambda\in\fa_{0,\C}^*$, restriction
\[
 \mcS'(X_0,\mcE_\lambda)^{[\fn]}\longrightarrow\mcS'(V,\mcE_\lambda|_V)^{[\fn]}
\]
is surjective.
\end{maintheorem}

Every tempered distribution on $V$ is the restriction of one on $X_0$ (see \cite[Theorem~5.4.1]{AG08}). The content of Theorem~\ref{thm:minimal-extension} is that the extension can be chosen $\IN$-finite.
An alternative, cohomological proof of
Theorem~\ref{thm:minimal-extension}
is given in \cite{XY26}.

Casselman--Hecht--Mili\v{c}i\'{c} developed the Bruhat filtration for Schwartz sections and tempered distributions and applied it to Whittaker functionals and the holomorphic continuation of Jacquet integrals \cite{CHM00}. Abe used the Bruhat filtration to study generalized Jacquet modules of parabolically induced representations \cite{Abe12}, i.e., the case $Q$ is minimal. His orbitwise description uses normal derivatives of polynomial densities, and the trivialization of equivariant bundles is given by a unipotent parametrization of each Bruhat cell. The analytic input in \cite[Proposition~4.4]{Abe12} includes the meromorphic continuation of Knapp--Stein intertwining operators. Deformation and pole removal then give the required extensions \cite[Lemmas~4.5--4.6]{Abe12}.

We first reduce Theorem~\ref{thm:minimal-extension} to the following statement by an induction over the finitely many $Q$-orbits on $X_0$.

\begin{maintheorem}\label{thm:supported-extension}
Let $O_w\subset X_0$ be a $Q$-orbit and $U_w\subset X_0$ a Nash open subset in which $O_w$ is closed. For every $\lambda\in\fa_{0,\C}^*$, restriction
\[
 \mcS'_{\overline{O_w}}(X_0,\mcE_\lambda)^{[\fn]}\longrightarrow\mcS'_{O_w}(U_w,\mcE_\lambda)^{[\fn]}
\]
is surjective. Here $\overline{O_w}$ is the (analytic) closure of $O_w$ in $X_0$.
\end{maintheorem}

\begin{proof}[Proof of Theorem~\ref{thm:minimal-extension} (assuming Theorem~\ref{thm:supported-extension})]
We induct on the finite number of $Q$-orbits in $V$ (Proposition~\ref{prop:flag-geometry}). The case $V=\varnothing$ is immediate. Otherwise choose an orbit $O_w\subset V$ of minimal dimension. Its analytic boundary $\overline{O_w}\setminus O_w$ is a union of lower-dimensional $Q$-orbits, so it is disjoint from the $Q$-stable set $V$. Thus $\overline{O_w}\cap V=O_w$, and $V'=V\setminus O_w$ is again a $Q$-stable Nash open subset. For $T\in\mcS'(V,\mcE_\lambda)^{[\fn]}$, induction gives $T_1\in\mcS'(X_0,\mcE_\lambda)^{[\fn]}$ extending $T|_{V'}$. The difference $T-T_1|_V$ is $\IN$-finite and supported in $O_w$, so Theorem~\ref{thm:supported-extension} with $U_w=V$ gives $T_2\in\mcS'_{\overline{O_w}}(X_0,\mcE_\lambda)^{[\fn]}$ extending it. Then $T_1+T_2$ extends $T$.
\end{proof}

\subsection{Method and organization}\label{ssec:organization}

We explain the idea behind the proof of Theorem~\ref{thm:supported-extension}. Distributions with support can be locally separated into normal derivatives and distributions along the support, where the derivative part is essentially global vector fields, of global nature.
To extend a distribution from a non-closed orbit $O_w$ to $X_0$, a natural first idea is to extend it to the closure $\overline{O_w}$ and then push it forward along the proper inclusion $\overline{O_w}\hookrightarrow X_0$. Formally, we would pair the distribution with test sections on $\overline{O_w}$, leaving boundary singularities. Two difficulties immediately arise: the orbit closure need not be smooth, so Schwartz functions and temperd distributions are not defined, and the boundary itself may be highly singular. Resolution of singularities addresses both geometric difficulties by replacing the closure with a smooth space whose boundary is a simple normal crossings divisor. The boundary is then locally a union of coordinate hyperplanes, suggesting a geometric approach to regularize the singular distributions by these coordinates, compared to Abe's approach, where continuation is in the inducing representation parameter $\lambda$, obtained through Knapp--Stein intertwining operators~\cite[Proposition~4.4]{Abe12}.

For non-minimal $Q$, however, $O_w$ need not be a single $N$-orbit. Rather, it has the form $O_w\simeq L\times^{R_w}F_w$, with $F_w=N/N_w$, and generalized $N$-invariance gives regularity only along the fibres. We therefore carry out this construction along the fibres: resolve the closure of $F_w$, then form the associated bundle to obtain a resolution of the orbit closure.

A further difficulty is that the resulting description by polynomial fibre densities uses specific $N$-invariant local frames. In these local frames, $\IN$-finite distributions on $O_w$ have polynomial density coefficients along the fibres, and distribution coefficients on the base (Proposition~\ref{prop:poly-densities}). A general tempered change of frame need not preserve polynomial coefficients along the fibre, and these orbit frames do not extend across the boundary. Comparing them with frames defined across the boundary introduces parameter-dependent tempered factors, which must be included in the analytic continuation step. The local extensions must also be glued compatibly with the bundle transitions.

Section~\ref{sec:conventions} fixes notation. Section~\ref{sec:normal} establishes the normal presentation for general $P$ and Casselman--Wallach coefficients. Sections~\ref{sec:fiber}--\ref{sec:bad} treat $P=P_0$, starting with the polynomial description of orbit distributions in Section~\ref{sec:fiber}. Section~\ref{sec:resolution} resolves the fibre closure equivariantly and compares the orbit frames with frames extending across the boundary. The resulting transition factors involve complex powers of algebraic functions and Nash coefficients. Section~\ref{sec:regularized} continues the local integrals meromorphically, using the resolution method of Atiyah~\cite{Ati70} and the idea of  rectilinearization of Bierstone--Milman~\cite{BM88} and Parusi\'nski~\cite{Par94}. Away from a locally finite union of parameter hyperplanes, Section~\ref{sec:zero-order} glues the regularized values into an $\fn$-equivariant tempered extension. Section~\ref{sec:bad} handles the remaining parameters by deformation and pole removal, following \cite[Lemmas~4.5--4.6]{Abe12}.
Section~\ref{sec:general-parabolics} proves coefficient lifting with support and deduces both assertions of Theorem~\ref{thm:general-extension}.

\section*{Acknowledgements}

This work is originally part of a research program initiated by the author's advisor. A complete proof was first obtained in June 2026, and the first draft of the manuscript was completed in July 2026.
The author is grateful to his advisor, Professor Jun Yu, for his guidance and valuable suggestions throughout the development of this work and the subsequent revisions of the manuscript. The author also thanks Kaidi Wu, Hongfeng Zhang, Cheng Chen for helpful discussions and suggestions.
The author acknowledges the assistance of GPT-6 Astra in improving the exposition and critically reviewing the proofs. The author takes full responsibility for the mathematical content and any remaining errors.

\section{Notation and conventions}\label{sec:conventions}

Algebraic groups and varieties are understood as schemes over $\R$ and are denoted by bold letters; their real points carry the analytic topology. By a variety we mean a geometrically integral separated scheme of finite type over $\R$. A variety with a smooth real point has Zariski-dense real points. Unless stated otherwise, closures of algebraic subvarieties are scheme-theoretic, whereas closures of subsets of real manifolds are analytic. When explicitly taking the Zariski closure of a subset of $\mathbb A^N_\R(\R)$, we mean the reduced closed subscheme defined by its ideal of vanishing polynomials.

\subsection{Groups}\label{ssec:groups}

Let $\G$ be a connected reductive algebraic group over $\R$, with Lie algebra $\fg$, and put $\fg_\C=\fg\otimes_\R\C$. Fix a finite-index subgroup $G\subset\G(\R)$ containing $\G(\R)^\circ$ and an algebraic Cartan involution $\tTheta$ preserving $G$. Write $\Theta=\tTheta|_G$ and $K=G^\Theta$. Parabolic subgroups of $G$ are understood to be $P=G\cap\pP(\R)$ for parabolic $\R$-subgroups $\pP\subset\G$.

Fix a maximal $\R$-split torus $\sS$ on which $\tTheta$ acts by inversion and a minimal parabolic $\pP_0\supset\sS$. Put $A_0=\sS(\R)^\circ$ and $\fa_0=\Lie\sS$. Identify $\Phi=\Phi(\G,\sS)$ with the restricted root system, with positive system $\Phi^+$ and simple roots $\Delta$ determined by $\pP_0$. For $I\subset\Delta$, write
\[
 \pP_I=\lL_I\ltimes\U_I,\qquad P_I=L_IN_I=M_IA_IN_I
\]
for the standard parabolic with its $\tTheta$-stable Levi subgroup and Langlands decomposition, where $L_I=G\cap\lL_I(\R)$ and $N_I=\U_I(\R)$. Put $\Phi_I=\Phi(\lL_I,\sS)$. A superscript $-$ denotes the opposite parabolic or unipotent radical with respect to the chosen Levi subgroup.

Let $W=N_G(A_0)/Z_G(A_0)$ be the relative Weyl group, choose representatives $\dot w\in N_G(A_0)$, and put $W_I=\langle s_\alpha:\alpha\in I\rangle$. Write ${}_IW_J$ for the minimal-length representatives of $W_I\backslash W/W_J$, and ${}_IW={}_IW_\emptyset$.

\subsection{Flag varieties and the Bruhat decomposition}\label{ssec:bruhat}

Fix standard parabolics $\pP=\pP_J=\lL_P\U_P$ and $\Q=\pP_I=\lL\U$, with real groups $P=L_PU_P$ and $Q=LN=MAN$. Put $\X=\G/\pP$, $X=G/P$, $\fn=\Lie N$, and $\IN=\ker(\mcU(\fn_\C)\to\C)$.

\begin{proposition}[real flag varieties and Bruhat decomposition]\label{prop:flag-geometry}
The natural maps
\[
 G/P\longrightarrow\G(\R)/\pP(\R)\longrightarrow\X(\R)
\]
are Nash diffeomorphisms of connected manifolds. For $w\in{}_IW_J$, put $\oO_w=\Q\dot w\pP/\pP$. Then
\begin{equation}\label{eq:Bruhat}
 X=\bigsqcup_{w\in{}_IW_J}O_w,\qquad
 O_w=\oO_w(\R)=Q\dot wP/P,
\end{equation}
and $\overline{O_w}\setminus O_w$ is a union of lower-dimensional $Q$-orbits. The translated big cells $\dot vU_P^-P/P$, indexed by $v\in W/W_J$, cover $X$.
\end{proposition}
\begin{proof}
The relative Bruhat decomposition gives these assertions for $\G(\R)$; see \cite[\S3.2, Theorem~3.13, Corollary~3.15 and Proposition~3.16]{BT72}. The real Bruhat cells have connected closures containing the base point, so $\X(\R)$ is connected. Since $\G(\R)^\circ$ has open orbits, it acts transitively. Applying this argument to the minimal flag variety of $\lL$ gives $\Q(\R)=\Q(\R)^\circ Z_\G(\sS)(\R)$. The second factor fixes $\dot w\pP$, and $\Q(\R)^\circ\subset Q$, so the $\Q(\R)$-orbit equals the $Q$-orbit.
\end{proof}

For $x_w=\dot w\pP\in\X(\R)$, multiplication gives the root-subgroup factorization
\begin{equation}\label{eq:N-factorization}
 (\U\cap\dot w\U_P^-\dot w^{-1})\times
 (\U\cap\dot w\pP\dot w^{-1})\xrightarrow{\ \sim\ }\U
\end{equation}
as $\R$-varieties. The orbit map identifies the first factor with $\U\cdot x_w$, hence
\begin{equation}\label{eq:N-orbit-real-points}
 (\U\cdot x_w)(\R)=N\cdot x_w.
\end{equation}

\subsection{Bundles and tempered distributions}\label{ssec:bundles}

For a real reductive group $H$, let $\mathrm{CW}(H)$ denote its Casselman--Wallach category. Fix $F\in\mathrm{CW}(L_P)$, extend it trivially across $U_P$, and put
\[
 \mcE_F=G\times^P F,\qquad [gp,v]=[g,pv].
\]
This is the unnormalized homogeneous bundle on $X$, with its canonical tempered Fr\'echet structure \cite[Section~3.3]{CS21}. For a Nash open subset $Y\subset X$, write $\mcS(Y,\mcE_F)$ for its Schwartz sections and $\mcS'(Y,\mcE_F)=\mcS(Y,\mcE_F)'_b$ for the strong dual. In a tempered trivialization, for example when $G\to X$ is trivial on $Y$, we have
\[
 \mcS(Y,\mcE_F)\simeq\mcS(Y)\widehat\otimes F,
\]
where $\widehat\otimes$ denotes the completed projective tensor product. The test-section spaces are nuclear Fr\'echet spaces. Multiplication by a tempered function is continuous on them. A short exact sequence of Casselman--Wallach coefficients induces a strict short exact sequence of test-section spaces.

For $P=P_0$, a finite-dimensional representation $\sigma$ of $M_0$ and $\lambda\in\fa_{0,\C}^*$, write $V_{\sigma,\lambda}$ for $V_\sigma$ with the $P_0$-action
\[
 \rho_{\sigma,\lambda}(man)v=e^{\lambda(\log a)}\sigma(m)v,
 \qquad \mcE_\lambda=G\times^{P_0}V_{\sigma,\lambda}.
\]
Since $M_0$ is compact, $\sigma$ is Nash \cite[Theorem~1.1(2)]{Sun15}.

The differential of the $G$-action makes $\mcU(\fg_\C)$ act by differential operators on sections of $\mcE_F$ over any open subset $Y$, continuously on $\mcS(Y,\mcE_F)$, and hence, through the transposes below, on $\mcS'(Y,\mcE_F)$. Restriction to smaller open subsets commutes with this action. We fix the transpose convention
\begin{equation}\label{eq:transpose}
 \xi^t=-\xi\quad(\xi\in\fg_\C),\qquad (D_1D_2)^t=D_2^tD_1^t,\qquad \langle DT,\phi\rangle=\langle T,D^t\phi\rangle.
\end{equation}
On scalar functions, $(\xi\phi)(x)=\frac{\dd}{\dd t}\big|_{t=0}\phi(\exp(-t\xi)x)$ for $\xi\in\fg$.

For a $\mcU(\fn_\C)$-module $E$ and $k\ge0$ we put
\[
 E^{\IN^k}=\{e\in E:\IN^ke=0\},\qquad E^{[\fn]}=\bigcup_{k\ge1}E^{\IN^k}.
\]
Support conditions are indicated by a subscript: $\mcS'_Z(Y,\mcE_F)$ consists of the elements of $\mcS'(Y,\mcE_F)$ supported in the closed subset $Z\subset Y$.

\begin{lemma}[independence of the neighbourhood]\label{lem:supported-nbhd}
Let $U'\subset U$ be Nash open subsets of $X$ with $O_w\subset U'$ and $O_w$ closed in $U$. Then restriction
\[
 \mcS'_{O_w}(U,\mcE_F)\longrightarrow\mcS'_{O_w}(U',\mcE_F)
\]
is a $\mcU(\fg_\C)$-module isomorphism.
\end{lemma}
\begin{proof}
A tempered partition of unity for $\{U',U\setminus O_w\}$ \cite[Theorem~5.2.1]{AG08} provides a tempered function $a$ on $U$, equal to one near $O_w$ and supported in $U'$, such that $\phi\mapsto(a\phi)|_{U'}$ is continuous on Schwartz sections. The inverse of restriction is $\langle T_U,\phi\rangle=\langle T_{U'},(a\phi)|_{U'}\rangle$.
\end{proof}

We write $\mcS'_{O_w}(X,\mcE_F)$ for the distributions supported on $O_w$ in $X\setminus(\overline{O_w}\setminus O_w)$. The lemma identifies this space canonically with $\mcS'_{O_w}(U,\mcE_F)$ for any Nash open subset $U\subset X$ in which $O_w$ is closed. For a closed embedding $i:O_w\hookrightarrow U$, the order-zero pushforward is defined by
\[
 \langle i_*u,\phi\rangle=\langle u,i^*\phi\rangle\qquad(u\in\mcS'(O_w,\mcE_F)).
\]

The following estimate follows from Kolchin's theorem (see \cite[Lemma~5.3 and Proposition~7.5]{Sun15}).

\begin{lemma}\label{lem:commutator}
\begin{enumerate}[label=\textup{(\alph*)}]
\item Let $E_1$ and $E_2$ be $\fn_\C$-modules, and let $\fn_\C$ act on $\operatorname{Hom}_\C(E_1,E_2)$ by $\xi\cdot A=\xi A-A\xi$. If $\IN^aA=0$ and $\IN^ke=0$, where $a,k\ge1$, then
\[
 \IN^{a+k-1}(Ae)=0.
\]
\item Every finite-dimensional rational representation of $N$ (that is, of $\U$) is annihilated by a power of $\IN$. In particular, for each $r\ge0$ there is $b_r\ge1$ such that $\IN^{b_r}$ annihilates $\mcU_{\le r}(\fg_\C)$ under the adjoint action, and then, for every $\fg_\C$-module $E$, $D\in\mcU_{\le r}(\fg_\C)$ and $k\ge1$,
\[
 \IN^ku=0\quad\Longrightarrow\quad\IN^{k+b_r-1}(Du)=0.
\]
\end{enumerate}
\end{lemma}

\section{Normal presentation along an orbit}\label{sec:normal}

Fix a $Q$-orbit $O_w\subset X=G/P$ and $F\in\mathrm{CW}(L_P)$. We construct local normal decompositions and a global normal presentation for distributions supported on $O_w$, and deduce a reduction of supported extension to normal order zero.

\subsection{The normal-jet filtration}\label{ssec:borel}

The following form of Borel's lemma is \cite[Proposition~8.2 and Proposition~8.3]{CS21}.

\begin{lemma}[Borel lemma]\label{lem:Borel}
Let a linear Nash group $H$ act on a Nash manifold $Y$, let $Z\subset Y$ be a closed $H$-stable Nash submanifold, and let $\mcE$ be an $H$-equivariant tempered vector bundle with nuclear Fr\'echet fibres on $Y$. Let
\[
 J^\infty_Z(Y,\mcE)=\mcS(Y,\mcE)/\mcS(Y\setminus Z,\mcE|_{Y\setminus Z})
\]
be the space of formal normal jets along $Z$, and put
\[
\begin{aligned}
 K_r&=\{\phi\in\mcS(Y,\mcE):\text{all derivatives of }\phi\text{ of order}<r\text{ vanish along }Z\},\\
 K_\infty&=\bigcap_{r\ge0}K_r.
\end{aligned}
\]
Derivatives are taken in local tempered trivializations of $\mcE$. Then
\[
 J^\infty_Z(Y,\mcE)=K_0/K_\infty\supset K_1/K_\infty\supset\cdots
\]
is a complete decreasing filtration by closed $H$-submodules, and
\[
 J^\infty_Z(Y,\mcE)\xrightarrow{\ \sim\ }\varprojlim_r\mcS(Y,\mcE)/K_r.
\]
Consequently the spaces $F_r\mcS'_Z(Y,\mcE)=\{u\in\mcS'(Y,\mcE):u|_{K_{r+1}}=0\}$ of distributions of normal order at most $r$ form an increasing filtration
\[
 F_0\mcS'_Z(Y,\mcE)\subset F_1\mcS'_Z(Y,\mcE)\subset\cdots,\qquad \mcS'_Z(Y,\mcE)=\bigcup_{r\ge0}F_r\mcS'_Z(Y,\mcE).
\]
\end{lemma}

Put $\mathcal N_{Z/Y}=(TY|_Z/TZ)_\C$ and $F_{-1}=0$. The normal Taylor map gives the canonical identification
\begin{equation}\label{eq:normal-symbol}
 \frac{F_r\mcS'_Z(Y,\mcE)}{F_{r-1}\mcS'_Z(Y,\mcE)}
 \simeq\mcS'\bigl(Z,\Sym^r\mathcal N_{Z/Y}^\vee\otimes\mcE|_Z\bigr),
\end{equation}
see \cite[Proposition~8.3]{CS21} and, for Nash bundles, \cite[Lemma~5.5.3 and Corollary~5.5.4]{AG08}. We write $\sigma_r$ for the resulting normal-symbol map. It is equivariant for the natural actions on the bundles.

The following is a product form of Schwartz's structure theorem for distributions supported on a submanifold, in the Nash setting; cf.\ the proof of \cite[Lemma~3.1]{SZ12}.

\begin{lemma}[local normal form]\label{lem:local-normal-form}
Let $Z$ be a Nash manifold, $Y=Z\times\R^d$ with normal coordinate $t\in\R^d$, and $E$ a trivial tempered bundle on $Y$ with nuclear Fr\'echet fibre. Every $T\in F_r\mcS'_Z(Y,E)$ has a unique expression
\[
 T=\sum_{|\alpha|\le r}\partial_t^\alpha T_\alpha,\qquad T_\alpha\in\mcS'(Z,E|_Z),
\]
and the coefficient maps $T\mapsto T_\alpha$ are continuous.
\end{lemma}

For a fixed $\chi\in C_c^\infty(\R^d)$ equal to $1$ near $0$, the coefficients are given by
\[
 \langle T_\alpha,\varphi\rangle
 =\frac{(-1)^{|\alpha|}}{\alpha!}
   \langle T,\chi(t)t^\alpha\varphi(z)\rangle,
 \qquad \varphi\in\mcS(Z,E|_Z).
\]

\subsection{The geometry of an orbit}\label{ssec:orbit-geometry}

We construct a neighbourhood $\mcU_w$ of $O_w$ covered by charts separating tangential and normal coordinates.

We take $w\in{}_IW_J$ (Section~\ref{ssec:bruhat}). Put
\[
\begin{gathered}
 x_w=\dot wP/P,\qquad H_w=Q\cap\dot wP\dot w^{-1}=R_w\ltimes N_w,\\
 R_w=L\cap\dot wP\dot w^{-1},\qquad N_w=N\cap\dot wP\dot w^{-1},
\end{gathered}
\]
with algebraic $\rR_w=\lL\cap\dot w\pP\dot w^{-1}$ and $\U_w=\U\cap\dot w\pP\dot w^{-1}$. The group $R_w$ is a standard parabolic of $L$. Then
\[
 O_w\cong Q/H_w\cong L\times^{R_w}F_w,\qquad F_w:=N\cdot x_w\cong N/N_w,
\]
where $R_w$ acts on $F_w$ by conjugation. Thus $O_w$ fibres over the real flag variety $L/R_w$ with fibre $F_w$. We write $\pi_w:O_w\to L/R_w$ for this map.

Let $P^-=L_PU_P^-$ and $Q^-=LN^-$ be the opposite parabolic subgroups. Put $\mathfrak u_P^-=\Lie U_P^-$, $\fp_P=\Lie P$, $\fn^-=\Lie N^-$, and
\[
 \fc_w=\Ad(\dot w)\mathfrak u_P^-\cap\fn^-.
\]
The relative root-space decomposition gives
\begin{equation}\label{eq:q-splitting}
 \fq=(\fq\cap\Ad(\dot w)\mathfrak u_P^-)\oplus(\fq\cap\Ad(\dot w)\fp_P),\qquad \Ad(\dot w)\mathfrak u_P^-=(\Ad(\dot w)\mathfrak u_P^-\cap\fq)\oplus\fc_w.
\end{equation}

\begin{lemma}\label{lem:tangent-normal}
The tangent space $T_{x_w}O_w$ and the normal space $N_{x_w}(O_w|X)$ are canonically identified with $\fq\cap\Ad(\dot w)\mathfrak u_P^-$ and $\fc_w$, respectively.
\end{lemma}

Let $\{B_\gamma\}_{\gamma\in W_L/W_{R_w}}$ be the cover of $L/R_w$ by translated big cells given by Proposition~\ref{prop:flag-geometry} for $L$, with fixed representatives $\dot\gamma\in L$. Write
\[
\begin{gathered}
 \B_\gamma=\dot\gamma(\lL\cap\dot w\U_P^-\dot w^{-1})\rR_w/\rR_w,\\
 B_\gamma=\B_\gamma(\R)=\dot\gamma(L\cap\dot wU_P^-\dot w^{-1})R_w/R_w.
\end{gathered}
\]
The algebraic sections
\[
 s_\gamma:\B_\gamma\to\lL,\qquad\dot\gamma l\rR_w\mapsto\dot\gamma l\qquad(l\in\lL\cap\dot w\U_P^-\dot w^{-1})
\]
take real points to $L$.
The sets $\Omega_\gamma=\pi_w^{-1}(B_\gamma)$ form an open cover of $O_w$, and
\begin{equation}\label{eq:Omega-product}
 \Omega_\gamma\xrightarrow{\ \sim\ }(L\cap\dot wU_P^-\dot w^{-1})\times F_w,\qquad[\dot\gamma nr,f]\mapsto(n,\Ad(r)f)\qquad(r\in R_w,\ f\in F_w).
\end{equation}

Next we extend each $\Omega_\gamma$ to an (affine) open subset $V_\gamma$ of $X$, in which we can separate tangential and normal data.

\begin{lemma}[normal charts]\label{lem:normal-charts}
Put
\[
 V_\gamma=\dot\gamma\dot wU_P^-P/P=\Ad(\dot\gamma\dot w)U_P^-\cdot\dot\gamma\dot wP/P,
\]
an open subset of $X$, and
\[
 \fc_{w,\gamma}=\Ad(\dot\gamma)\fc_w,\qquad \fl_{w,\gamma}=\Ad(\dot\gamma)(\Ad(\dot w)\mathfrak u_P^-\cap\fl),\qquad \fu_{w,\gamma}=\Ad(\dot\gamma)(\Ad(\dot w)\mathfrak u_P^-\cap\fn).
\]
Then
\begin{equation}\label{eq:normal-chart}
 \Phi_{w,\gamma}:\fc_{w,\gamma}\oplus\fl_{w,\gamma}\oplus\fu_{w,\gamma}\to V_\gamma,\qquad(C,Y)\mapsto\exp(C+Y)\,\dot\gamma\dot wP/P,
\end{equation}
where $\exp$ is the exponential map of the unipotent group $\Ad(\dot\gamma\dot w)U_P^-$, $C\in\fc_{w,\gamma}$ and $Y\in\fl_{w,\gamma}\oplus\fu_{w,\gamma}$, is an isomorphism of affine varieties, and
\[
 V_\gamma\cap O_w=\Omega_\gamma=\Phi_{w,\gamma}(\{0\}\oplus\fl_{w,\gamma}\oplus\fu_{w,\gamma}).
\]
In particular, $O_w$ is closed in $\bigcup_\gamma V_\gamma$, and in each $V_\gamma$ it is cut out by the $\fc_{w,\gamma}$-coordinates.
\end{lemma}

\begin{proof}
By \eqref{eq:q-splitting}, $\fc_{w,\gamma}\oplus\fl_{w,\gamma}\oplus\fu_{w,\gamma}=\Ad(\dot\gamma\dot w)\mathfrak u_P^-$. The first assertion holds because the exponential map of a unipotent group, and the orbit map of $\Ad(\dot\gamma\dot w)U_P^-$ through $\dot\gamma\dot wP$, are isomorphisms.

The rest relies on two factorizations:
\begin{align}
 Q\cap(\dot wU_P^-\dot w^{-1})(\dot wP\dot w^{-1})&=(Q\cap\dot wU_P^-\dot w^{-1})(Q\cap\dot wP\dot w^{-1}),\label{eq:claimI}\\
 N&=(N\cap\dot wU_P^-\dot w^{-1})(N\cap\dot wP\dot w^{-1}).\label{eq:claimII}
\end{align}
For \eqref{eq:claimI}: multiplication $m:\dot w\U_P^-\dot w^{-1}\times\dot w\pP\dot w^{-1}\to(\dot w\U_P^-\dot w^{-1})(\dot w\pP\dot w^{-1})$ is an isomorphism onto an open subset of $\G$. Hence $\Q\cap(\dot w\U_P^-\dot w^{-1})(\dot w\pP\dot w^{-1})$ is smooth and geometrically integral, of dimension $\dim\fq$. It contains the image $m\bigl((\Q\cap\dot w\U_P^-\dot w^{-1})\times(\Q\cap\dot w\pP\dot w^{-1})\bigr)$ as a closed subscheme of dimension $\dim(\fq\cap\Ad(\dot w)\mathfrak u_P^-)+\dim(\fq\cap\Ad(\dot w)\fp_P)$, which is $\dim\fq$ by \eqref{eq:q-splitting}. So the two coincide, and taking real points and intersecting with $G$ gives \eqref{eq:claimI}.
Factorization \eqref{eq:claimII} is the real-point version of \eqref{eq:N-factorization}.

Now let $\dot\gamma n\dot wP\in V_\gamma\cap O_w$ with $n\in\dot wU_P^-\dot w^{-1}$. Since $\dot\gamma\in L\subset Q$, we have $n\in Q\cdot\dot wP\dot w^{-1}$; write $n=qp^{-1}$ with $q\in Q$ and $p\in\dot wP\dot w^{-1}$. By \eqref{eq:claimI},
\[
 q=np\in Q\cap(\dot wU_P^-\dot w^{-1})(\dot wP\dot w^{-1})=(Q\cap\dot wU_P^-\dot w^{-1})(Q\cap\dot wP\dot w^{-1}),
\]
say $q=q_1q_2$. Then $n=q_1\cdot(q_2p^{-1})$ is a second factorization of $n$ in the open cell, besides $n=n\cdot1$. Uniqueness gives $n=q_1\in Q\cap\dot wU_P^-\dot w^{-1}$, and hence $\dot\gamma n\dot wP=(\Ad(\dot\gamma)n)\,\dot\gamma\dot wP\in\Phi_{w,\gamma}(\{0\}\oplus\fl_{w,\gamma}\oplus\fu_{w,\gamma})$. Conversely, $\dot\gamma n\dot wP\in V_\gamma\cap O_w$ for every $n\in Q\cap\dot wU_P^-\dot w^{-1}$. Therefore
\[
 V_\gamma\cap O_w=\Phi_{w,\gamma}(\{0\}\oplus\fl_{w,\gamma}\oplus\fu_{w,\gamma})=\dot\gamma(Q\cap\dot wU_P^-\dot w^{-1})\dot wP/P.
\]

It remains to identify this set with $\Omega_\gamma$. Under $O_w\cong Q/H_w$, $q\dot wP\mapsto qH_w$, the map $\pi_w$ is the projection to the Levi factor, and $\pr_L(Q\cap\dot wU_P^-\dot w^{-1})=L\cap\dot wU_P^-\dot w^{-1}$. Hence $\pi_w(\dot\gamma n\dot wP)=\dot\gamma\pr_L(n)R_w\in B_\gamma$ for $n\in Q\cap\dot wU_P^-\dot w^{-1}$, that is, $V_\gamma\cap O_w\subset\Omega_\gamma$. Conversely, let $lu\dot wP\in\Omega_\gamma$ with $l\in L$ and $u\in N$. Since $\pi_w(lu\dot wP)=lR_w\in B_\gamma$, we can write $l=\dot\gamma nr$ with $n\in L\cap\dot wU_P^-\dot w^{-1}$ and $r\in R_w$, and then
\[
 lu\dot wP=\dot\gamma n(rur^{-1})\dot w(\dot w^{-1}r\dot w)P=\dot\gamma n(rur^{-1})\dot wP.
\]
By \eqref{eq:claimII}, $rur^{-1}=u^wu_w$ with $u^w\in N\cap\dot wU_P^-\dot w^{-1}$ and $u_w\in N\cap\dot wP\dot w^{-1}$. Hence $lu\dot wP=\dot\gamma(nu^w)\dot wP$ with $nu^w\in Q\cap\dot wU_P^-\dot w^{-1}$, which proves $\Omega_\gamma\subset V_\gamma\cap O_w$.
\end{proof}

\begin{remark}\label{rem:cw}
The algebraic factorization underlying \eqref{eq:claimII} shows that multiplication $(\U\cap\dot w\U_P^-\dot w^{-1})\times\U_w\to\U$ is an isomorphism. Let $\F_w=\U\dot w\pP/\pP\cong\U/\U_w$. We obtain an isomorphism
\[
 c_w:\F_w\xrightarrow{\ \sim\ }\U\cap\dot w\U_P^-\dot w^{-1},\qquad c_w(y)\dot w\pP=y.
\]
By \eqref{eq:N-orbit-real-points}, $\F_w(\R)=N\cdot x_w=F_w$.
\end{remark}

\subsection{Normal presentation and reduction to order zero}\label{ssec:order-zero}

We first state the local decomposition and global presentation, then deduce the reduction to normal order zero. The local and global results are proved in Sections~\ref{ssec:local-presentation} and~\ref{ssec:global-presentation}, respectively.

\begin{definition}\label{def:Uw}
Put $\mcU_w=\bigcup_{\gamma\in W_L/W_{R_w}}V_\gamma$; by Lemma~\ref{lem:normal-charts}, $O_w$ is closed in $\mcU_w$. Put
\[
\begin{aligned}
 F_r\mcS'_{O_w}(\mcU_w,\mcE_F)^{[\fn]}&=\bigcup_{k\ge0}F_r\mcS'_{O_w}(\mcU_w,\mcE_F)^{\IN^{k+1}},\\
 \mcS'(O_w,\mcE_F)^{[\fn]}&=\bigcup_{k\ge0}\mcS'(O_w,\mcE_F)^{\IN^{k+1}}.
\end{aligned}
\]
\end{definition}

\begin{proposition}[local normal decomposition]\label{prop:local-presentation}
For every $r\ge0$, the map
\[
 \Theta_{\gamma,r}:\mcU_{\le r}((\fc_{w,\gamma})_\C)\otimes\mcS'(\Omega_\gamma,\mcE_F)\longrightarrow F_r\mcS'_{\Omega_\gamma}(V_\gamma,\mcE_F),\qquad D\otimes u\mapsto Du,
\]
where $\langle Du,\phi\rangle=\langle u,(D^t\phi)|_{\Omega_\gamma}\rangle$, is a topological isomorphism (LHS carries the finite direct-sum topology). Fix an ordered basis of $\fc_{w,\gamma}$, with associated PBW basis $\{D_i\}$ of $\mcU_{\le r}((\fc_{w,\gamma})_\C)$, and write $T=\sum_iD_iu_i$. If $\supp T\subset C_\gamma$ for a closed subset $C_\gamma\subset\Omega_\gamma$, then $\supp u_i\subset C_\gamma$ for every $i$.
\end{proposition}

\begin{lemma}[preservation of $\IN$-finiteness]\label{lem:IN-bounds}
For every $r\ge0$, the map $\Theta_{\gamma,r}$ of Proposition~\ref{prop:local-presentation} restricts to an isomorphism
\[
 \mcU_{\le r}((\fc_{w,\gamma})_\C)\otimes\mcS'(\Omega_\gamma,\mcE_F)^{[\fn]}\xrightarrow{\ \sim\ }F_r\mcS'_{\Omega_\gamma}(V_\gamma,\mcE_F)^{[\fn]}.
\]
Moreover, for fixed $r\ge0$ and $m\ge1$, a common power of $\IN$, independent of $T$ and $F$, annihilates all the coefficients $u_i$ in a fixed PBW expansion $T=\sum_iD_iu_i$ whenever $\IN^mT=0$.
\end{lemma}

The local decompositions yield the following global presentation.

\begin{theorem}[global normal presentation]\label{thm:global-presentation}
The map
\[
\begin{gathered}
 \Pi_r:\bigoplus_\gamma\mcU_{\le r}((\fc_{w,\gamma})_\C)\otimes\mcS'(O_w,\mcE_F)^{[\fn]}\longrightarrow F_r\mcS'_{O_w}(\mcU_w,\mcE_F)^{[\fn]},\\
 D\otimes u\longmapsto\bigl[\phi\mapsto\langle u,(D^t\phi)|_{O_w}\rangle\bigr]\qquad(\phi\in\mcS(\mcU_w,\mcE_F)),
\end{gathered}
\]
is surjective. Moreover, the maps $\Pi_r$ admit compatible linear sections $s_r$ satisfying
\[
 \Pi_rs_r=\id,\qquad s_{r+1}|_{F_r}=s_r.
\]
The coefficients in the $\gamma$-summand of $s_r(T)$ are supported in $C_\gamma=\pi_w^{-1}(K_\gamma)$, for fixed compact subsets $K_\gamma\subset B_\gamma$ independent of $F$.
\end{theorem}

The global presentation reduces supported extension to normal order zero. Let $r^X_w:\mcS'(X,\mcE_F)\to\mcS'(\mcU_w,\mcE_F)$ be the restriction and $i:O_w\hookrightarrow\mcU_w$ the closed embedding.

\begin{proposition}[reduction to order zero]\label{prop:order-zero}
Let
\[
 \E_0:\mcS'(O_w,\mcE_F)^{[\fn]}\longrightarrow\mcS'_{\overline{O_w}}(X,\mcE_F)^{[\fn]}
\]
be a linear map with $r^X_w\circ\E_0=i_*$. For $r\ge0$ define
\[
 \widetilde\E_r:\bigoplus_\gamma\mcU_{\le r}((\fc_{w,\gamma})_\C)\otimes\mcS'(O_w,\mcE_F)^{[\fn]}\to\mcS'_{\overline{O_w}}(X,\mcE_F)^{[\fn]},\qquad D\otimes u\mapsto D\cdot\E_0(u),
\]
and $\E_r=\widetilde\E_r\circ s_r$ on $F_r\mcS'_{O_w}(\mcU_w,\mcE_F)^{[\fn]}$. Then $\E_{r+1}$ extends $\E_r$, and their union
\[
 \E_{w,F}:\mcS'_{O_w}(\mcU_w,\mcE_F)^{[\fn]}\longrightarrow\mcS'_{\overline{O_w}}(X,\mcE_F)^{[\fn]}
\]
satisfies $r^X_w\circ\E_{w,F}=\id$.
\end{proposition}

\begin{proof}
The image of $\widetilde\E_r$ is $\IN$-finite by Lemma~\ref{lem:commutator}(b), and it is supported in $\overline{O_w}$. Restriction commutes with $D$, so $(D\E_0(u))|_{\mcU_w}=Di_*u$; summing over the presentation $\Pi_rs_r(T)=T$ gives $r^X_w\circ\E_r=\id$ on $F_r\mcS'_{O_w}(\mcU_w,\mcE_F)^{[\fn]}$. Since $s_{r+1}|_{F_r}=s_r$, the maps $\E_r$ are compatible, and every element of $\mcS'_{O_w}(\mcU_w,\mcE_F)$ has finite normal order (Lemma~\ref{lem:Borel}). Thus $\E_{w,F}$ is the required section. The same conclusion holds in any Nash open neighbourhood in which $O_w$ is closed, by Lemma~\ref{lem:supported-nbhd}.
\end{proof}

\subsection{Proof of the local normal decomposition}\label{ssec:local-presentation}

We prove Proposition~\ref{prop:local-presentation} and Lemma~\ref{lem:IN-bounds}. First we compare Lie algebra derivatives with coordinate normal derivatives.

\begin{lemma}\label{lem:bundle-trivial}
The map
\[
 \Psi_\gamma:V_\gamma\times F\to\mcE_F|_{V_\gamma},\qquad(\dot\gamma\dot wnP,v)\mapsto[\dot\gamma\dot wn,v]\qquad(n\in U_P^-,\ v\in F),
\]
is a tempered trivialization of $\mcE_F$ over $V_\gamma$. Under the resulting isomorphism $\mcS(V_\gamma,\mcE_F|_{V_\gamma})\cong\mcS(V_\gamma)\widehat\otimes F$, each $\xi\in\fc_{w,\gamma}$ acts by its natural action on the first factor and trivially on the second.
\end{lemma}

\begin{proof}
The trivialization comes from a Nash section, so it is tempered by \cite[Section~3.3]{CS21}. For $x=\dot\gamma\dot wnP$ the chosen lift is $g_x=\dot\gamma\dot wn$. If $\xi\in\fc_{w,\gamma}\subset\Ad(\dot\gamma\dot w)\mathfrak u_P^-$, then $\exp(-t\xi)g_x\in\dot\gamma\dot wU_P^-$, so it is the chosen lift of $\exp(-t\xi)x$. Hence the induced action differentiates along the flow $\exp(-t\xi)$ and does not act on $F$.
\end{proof}

Fix $\gamma$ and put $\fc=\Ad(\dot w)\mathfrak u_P^-$. Conjugate the coordinates of \eqref{eq:normal-chart} by $\Ad(\dot\gamma)^{-1}$ and write them as $(C,Y)\in\fc_w\oplus(\fc\cap\fq)$, so that $\Omega_\gamma=\{C=0\}$. 
For $Z=C+Y$ and function $f$ on $V_\gamma$, put
\[
 x_Z=\dot\gamma\exp(Z)\dot wP,\qquad
\widetilde f(Z)=\widetilde f(C,Y)=f(x_Z).
\]
We first choose a normal frame and compare its PBW derivatives with coordinate normal derivatives.

\begin{lemma}[triangular normal frame]\label{lem:triangular-normal-frame}
There is an ordered root-space basis $e_1,\dots,e_s$ of $\fc$ whose first $d=\dim\fc_w$ vectors span $\fc_w=\fc\cap\fn^-$, such that, for $E_i=\Ad(\dot\gamma)e_i$,
\begin{equation}\label{eq:first-order}
 \widetilde{E_if}=-\Bigl(\partial_{e_i}+\sum_{j>i}a_{ij}(C,Y)\,\partial_{e_j}\Bigr)\widetilde f,
 \qquad 1\le i\le d,
\end{equation}
where $a_{ij}$ are polynomials. All $E_i$, $1\le i\le s$, have polynomial coefficients. The normal classes of $E_1,\dots,E_d$ along $\Omega_\gamma$ form an algebraic frame.
\end{lemma}

\begin{proof}
Put $h_w(\alpha)=\Ht(-w^{-1}\alpha)>0$ for $\alpha\in w(\Phi^-\setminus\Phi_J)$. Choose $e_i\in\fg_{\alpha_i}$ such that
\[
\begin{gathered}
 \fc_w=\Span\{e_1,\dots,e_d\},\qquad
 \fc\cap\fq=\Span\{e_{d+1},\dots,e_s\},\\
 h_w(\alpha_i)\le h_w(\alpha_j)
 \quad\text{if }1\le i<j\le d\text{ or }d<i<j\le s.
\end{gathered}
\]
We have $(\ad Z)^\nu e_i\in\Span\{e_j:j>i\}
 \quad(Z\in\fc,\ \nu\ge1,\ 1\le i\le d)$.

For $Z=C+Y$, direct computation gives
\begin{equation}\label{eq:BCH}
\begin{aligned}
 (E_if)(x_Z)
 &=\frac{\dd}{\dd t}\Big|_{t=0}
 f\bigl(\exp(-tE_i)\dot\gamma\exp(Z)\dot wP\bigr)\\
  &=\frac{\dd}{\dd t}\Big|_{t=0}
 f\bigl(\dot\gamma\exp(-te_i)\exp(Z)\dot wP\bigr)\\
 &=\frac{\dd}{\dd t}\Big|_{t=0}
 \widetilde f\bigl(\log(\exp(-te_i)\exp(Z))\bigr)\\
  &=\dd\widetilde f_Z\Bigl(
 \frac{\dd}{\dd t}\Big|_{t=0}\log(\exp(-te_i)\exp(Z))\Bigr)\\
 &=-\dd\widetilde f_Z\Bigl(
 \sum_{\nu\ge0}\frac{B_\nu}{\nu!}(\ad Z)^\nu e_i\Bigr),
\end{aligned}
\end{equation}
where $z/(e^z-1)=\sum_{\nu\ge0}B_\nu z^\nu/\nu!$. Nilpotence makes the sum finite, and each $(\ad Z)^\nu e_i$ has polynomial coordinates in $Z$. Thus every $E_i$ has polynomial coefficients on the whole chart.

Let $\overline E_i$ and $\overline{\partial_{C_i}}$ denote the normal classes along $\Omega_\gamma$. Setting $C=0$ and discarding tangential derivatives gives
\[
 \overline E_i=-\overline{\partial_{C_i}}
 -\sum_{i<j\le d}a_{ij}(0,Y)\overline{\partial_{C_j}},
 \qquad 1\le i\le d.
\]
The transition matrix between $\{\bar{E_i},1\leq i\leq d\}$ and coordinate derivatives is an upper triangular matrix with diagonal $-1$, which has polynomial inverse matrix.
\end{proof}

\begin{lemma}\label{lem:PBW-normal-jets}
Use the basis of Lemma~\ref{lem:triangular-normal-frame} and put $D_\alpha=E_1^{\alpha_1}\cdots E_d^{\alpha_d}$ for $\alpha\in\mathbb N^d$ (Expand $D_\alpha$ into sequence of basis ordered by $E_1<\cdots<E_d$. Then order multi-indices lexicographically). Then
\begin{equation}\label{eq:principal-term}
\begin{aligned}
 (D_\alpha^tf)|_{\Omega_\gamma}
 &=\partial_C^\alpha\widetilde f(0,Y)
 +\sum_{\substack{|\beta|=|\alpha|\\\beta>\alpha}}a_{\alpha\beta}(Y)\,\partial_C^\beta\widetilde f(0,Y)\\
 &\quad+\sum_{|\beta|<|\alpha|}A_{\alpha\beta}(Y,\partial_Y)\bigl(\partial_C^\beta\widetilde f(0,Y)\bigr),
\end{aligned}
\end{equation}
where the $a_{\alpha\beta}$ are polynomials and the $A_{\alpha\beta}$ are tangential differential operators with polynomial coefficients.
\end{lemma}

\begin{proof}
The principal symbol of a product is the product of the principal symbols. By \eqref{eq:first-order} at $C=0$, its terms consisting only of normal derivatives give $\partial_C^\alpha$ or a later multi-index of the same degree. A term containing a tangential derivative, or arising when a coefficient is differentiated, has lower normal degree. All coefficients remain polynomial, giving \eqref{eq:principal-term}.
\end{proof}

\begin{proof}[Proof of Proposition~\ref{prop:local-presentation}]
Use the trivialization of Lemma~\ref{lem:bundle-trivial}. The scalar Taylor maps and polynomial differential operators above act continuously after completed tensor product with $F$.

By Lemma~\ref{lem:local-normal-form}, the map $T\mapsto(v_\alpha)_{|\alpha|\le r}$ defined by
\[
 \langle T,f\rangle=\sum_{|\alpha|\le r}\langle v_\alpha,\partial_C^\alpha\widetilde f(0,\cdot)\rangle
\]
is a topological isomorphism. Lemma~\ref{lem:PBW-normal-jets} gives the inverse
\[
 \Theta_{\gamma,r}^{-1}:T\longmapsto\sum_{|\alpha|\le r}D_\alpha\otimes u_\alpha,
\]
where the $u_\beta$ are determined recursively, by decreasing degree and then by the order, from
\[
 u_\beta=v_\beta
 -\sum_{\substack{|\alpha|=|\beta|\\\alpha<\beta}}a_{\alpha\beta}u_\alpha
 -\sum_{|\alpha|>|\beta|}A_{\alpha\beta}'u_\alpha.
\]
Here $\langle A_{\alpha\beta}'u,\psi\rangle=\langle u,A_{\alpha\beta}\psi\rangle$. Polynomial-coefficient differential operators act continuously on Schwartz distributions, so both $\Theta_{\gamma,r}$ and this inverse are continuous.

The coefficient formula following Lemma~\ref{lem:local-normal-form} gives $\supp v_\alpha\subset C_\gamma$ if $\supp T\subset C_\gamma$. Since differential operators do not enlarge support, the same recursion gives $\supp u_\alpha\subset C_\gamma$.
\end{proof}

\begin{proof}[Proof of Lemma~\ref{lem:IN-bounds}]
We first extract coefficients from an arbitrary homogeneous normal symbol. Put $\mathcal N_\gamma=(TX|_{\Omega_\gamma}/T\Omega_\gamma)_\C$, and let
\[
 q\in\mcS'\bigl(\Omega_\gamma,\Sym^r\mathcal N_\gamma^\vee\otimes\mcE_F\bigr),\qquad \IN^mq=0.
\]
By Lemma~\ref{lem:triangular-normal-frame}, the transposed normal symbols of the PBW monomials $D_i$ with $\deg D_i=r$ form an algebraic frame of $\Sym^r\mathcal N_\gamma$. Let $\eta_i$ be the dual frame, with the Taylor normalization of Lemma~\ref{lem:PBW-normal-jets}, and define the coefficients $u_i\in\mcS'(\Omega_\gamma,\mcE_F)$ by
\[
 \langle u_i,\varphi\rangle
 =\langle q,\eta_i\otimes\varphi\rangle,
 \qquad\varphi\in\mcS(\Omega_\gamma,\mcE_F).
\]
Multiplication by $\eta_i$ is continuous on Schwartz sections, so these coefficient maps are continuous and preserve support. The $\eta_i$ are regular sections on the $\U$-stable affine algebraic chart underlying $\Omega_\gamma$. Since $\mathcal N_\gamma$ is an algebraic $\U$-equivariant bundle, they lie in a common finite-dimensional rational $\U$-module. By Lemma~\ref{lem:commutator}(b), there is $b\ge1$, depending only on $r$ and the fixed geometric data, such that $\IN^b\eta_i=0$ for all $i$.

For $X\in\fn_\C$, the tensor-product action gives
\[
\begin{aligned}
 \langle Xu_i,\varphi\rangle
 &=-\langle q,\eta_i\otimes X\varphi\rangle\\
 &=\langle Xq,\eta_i\otimes\varphi\rangle
   +\langle q,(X\eta_i)\otimes\varphi\rangle.
\end{aligned}
\]
Iterating this identity gives $\IN^{m+b-1}u_i=0$.

Now write $T=\sum_{\deg D_i\le r}D_iu_i$ with $\IN^mT=0$. By \eqref{eq:normal-symbol}, $\IN^m\sigma_r(T)=0$, and its coefficients above are precisely the $u_i$ with $\deg D_i=r$. Lemma~\ref{lem:commutator}(b) therefore makes the highest-degree terms $D_iu_i$ $\IN$-finite, and
\[
 T-\sum_{\deg D_i=r}D_iu_i
 \in F_{r-1}\mcS'_{\Omega_\gamma}(V_\gamma,\mcE_F)
\]
is $\IN$-finite. The bounds for the dual frames and for the adjoint action on $\mcU_{\le r}(\fg_\C)$ are independent of $T$ and $F$. Induction on $r$, starting with $F_{-1}=0$, proves the uniform assertion. Conversely, if all the $u_i$ are $\IN$-finite, Lemma~\ref{lem:commutator}(b) shows that each $D_iu_i$ is $\IN$-finite.
\end{proof}

\subsection{Proof of the global normal presentation}\label{ssec:global-presentation}

We localize the highest normal symbol on the base and then induct on normal order. Since $L/R_w$ is compact, choose a smooth partition of unity $\sum_\gamma\beta_\gamma=1$ with
\[
 \supp\beta_\gamma\subset K_\gamma\Subset B_\gamma,\qquad
 C_\gamma=\pi_w^{-1}(K_\gamma),
\]
where the $K_\gamma$ are fixed compact subsets. The functions $\beta_\gamma\circ\pi_w$ are tempered and $N$-invariant on $O_w$. These choices are independent of $r$ and $F$.

The following is an immediate consequence of Schwartz localization with an $N$-invariant cutoff.

\begin{lemma}[extension of supported coefficients]\label{lem:coefficient-extension}
Let $K\subset B_\gamma$ be compact and put $C=\pi_w^{-1}(K)$. Choose $\eta_\gamma\in C_c^\infty(B_\gamma)$ equal to one near $K$, extend it by zero to $L/R_w$, and put $\theta_\gamma=\eta_\gamma\circ\pi_w$. The formula
\[
 \langle\widetilde u,\phi\rangle
 =\langle u,(\theta_\gamma\phi)|_{\Omega_\gamma}\rangle,\qquad \phi\in\mcS(O_w,\mcE_F)
\]
defines a continuous $N$-equivariant extension map
\[
 \mcS'_C(\Omega_\gamma,\mcE_F)
 \longrightarrow\mcS'_C(O_w,\mcE_F),\qquad u\mapsto\widetilde u.
\]
The functions $\eta_\gamma$ and $\theta_\gamma$ can be fixed independently of $r$ and $F$.
\end{lemma}

\begin{proof}[Proof of Theorem~\ref{thm:global-presentation}]
Let $T\in F_r\mcS'_{O_w}(\mcU_w,\mcE_F)$ satisfy $\IN^kT=0$, and put $\mathcal N_w=(TX|_{O_w}/TO_w)_\C$. By \eqref{eq:normal-symbol}, its highest normal symbol lies in $\mcS'(O_w,\Sym^r\mathcal N_w^\vee\otimes\mcE_F)$. Set
\[
 q_\gamma=\bigl((\beta_\gamma\circ\pi_w)\sigma_r(T)\bigr)|_{\Omega_\gamma}.
\]
Since $N$ fixes the base, $\IN^kq_\gamma=0$. Using the dual frame $\eta_{\gamma i}$ from the proof of Lemma~\ref{lem:IN-bounds}, define $u_{\gamma i}\in\mcS'(\Omega_\gamma,\mcE_F)$, for $\deg D_{\gamma i}=r$, by
\[
 \langle u_{\gamma i},\varphi\rangle
 =\langle q_\gamma,\eta_{\gamma i}\otimes\varphi\rangle,
 \qquad \varphi\in\mcS(\Omega_\gamma,\mcE_F).
\]
Equivalently, if
\[
 T|_{V_\gamma}=\sum_{\deg D_{\gamma i}\le r}D_{\gamma i}v_{\gamma i}
\]
is the local PBW decomposition, then
\[
 u_{\gamma i}=(\beta_\gamma\circ\pi_w)v_{\gamma i},
 \qquad \deg D_{\gamma i}=r.
\]
The normal-symbol calculation in that proof shows that the $u_{\gamma i}$ are $\IN$-finite, and their support lies in $C_\gamma$. Extend them by Lemma~\ref{lem:coefficient-extension} and put
\[
 a_r(T)=\sum_{\gamma,\,\deg D_{\gamma i}=r}
 D_{\gamma i}\otimes\widetilde u_{\gamma i},
 \qquad T_r=\Pi_ra_r(T).
\]
Lemma~\ref{lem:commutator}(b) makes $T_r$ $\IN$-finite. Its normal symbol is the sum of the localized symbols, so
\[
 \sigma_r(T_r)=\sigma_r(T),\qquad
 T-T_r\in F_{r-1}\mcS'_{O_w}(\mcU_w,\mcE_F)^{[\fn]}.
\]
Induction on $r$ proves surjectivity. Compatible linear sections are obtained by the recursion below.
\end{proof}

\begin{remark}[the linear section $s_r$]\label{rem:section-sr}
Fix the base partition above, the cutoffs of Lemma~\ref{lem:coefficient-extension}, and the ordered PBW bases of Lemma~\ref{lem:triangular-normal-frame}, independently of $r$ and $F$. With $a_r$ as in the proof, define
\[
 s_{-1}=0,\qquad
 s_r(T)=a_r(T)+s_{r-1}\bigl(T-\Pi_ra_r(T)\bigr).
\]
Then $\Pi_rs_r=\id$. If $T\in F_{r-1}$, then $a_r(T)=0$, so $s_r|_{F_{r-1}}=s_{r-1}$. The coefficients in the $\gamma$-summand are supported in $C_\gamma$, and the construction is continuous on each $F_r\mcS'_{O_w}(\mcU_w,\mcE_F)^{\IN^k}$. Iterating the bounds in Lemmas~\ref{lem:IN-bounds} and~\ref{lem:commutator} gives a common exponent $K(w,r,k)\ge k$, independent of $T$ and $F$, annihilating every coefficient of $s_r(T)$.
\end{remark}

\section{\texorpdfstring{$\IN$}{I-N}-finite distributions on an orbit}\label{sec:fiber}

Throughout Sections~\ref{sec:fiber}--\ref{sec:bad}, we take $P=P_0$ and $F=V_{\sigma,\lambda}$, so that $\mcE_F=\mcE_\lambda$, and abbreviate $X_0$ to $X$.

This section describes normal order zero distributions $\mcS'(O_w,\mcE_\lambda)^{[\fn]}$. Under $N$-invariant frames over $\Omega_\gamma$ described below, an $\IN$-finite distribution is a distribution on the base tensored with a polynomial density along the fibre $F_w$ (Proposition~\ref{prop:poly-densities}).

\subsection{Orbit frames}\label{ssec:orbit-frames}

Fix a basis $e_1,\dots,e_d$ of $V_\sigma$, which is the underlying space of $V_{\sigma,\lambda}$ for every $\lambda$, and let $e^1,\dots,e^d$ be the dual basis. Put $\mathbf e=(e_1,\dots,e_d)$ and $\mathbf e^\vee=(e^1,\dots,e^d)$. The notation $[g,\mathbf e]$ denotes the row with entries $[g,e_j]$, and similarly for $[g,\mathbf e^\vee]$. Throughout, frames and dual frames are row vectors, coefficient vectors are columns, and $A^{-t}=(A^{-1})^t$ denotes the inverse transpose.

Using \eqref{eq:Omega-product}, we write points of $\Omega_\gamma$ as pairs $(b,f)\in B_\gamma\times F_w$, the pair $(b,uN_w)$ corresponding to $s_\gamma(b)u\dot wP_0$. The section $s_\gamma$ and the map $c_w$ of Remark~\ref{rem:cw} define the \emph{orbit lift}
\[
 g_\gamma:\B_\gamma\times\F_w\to\G,\qquad(b,f)\mapsto s_\gamma(b)c_w(f)\dot w.
\]
Define
\begin{equation}\label{eq:tau-frame}
 \tau_\gamma(b,f)=[g_\gamma(b,f),\mathbf e].
\end{equation}
For any representative $u\in N$ of $f=uN_w$, we also have $\tau_\gamma(b,f)=[s_\gamma(b)u\dot w,\mathbf e]$, since $\dot w^{-1}\U_w\dot w\subset\U_0$ acts trivially on $V_{\sigma,\lambda}$.
Write $\tau_\gamma=(\tau_{\gamma,1},\dots,\tau_{\gamma,d})$ and let $\tau_\gamma^\vee=(\tau_\gamma^1,\dots,\tau_\gamma^d)$ be the tempered dual frame, characterized by $\langle\tau_\gamma^i,\tau_{\gamma,j}\rangle=\delta^i_j$. Every $\phi\in\mcS(\Omega_\gamma,\mcE_\lambda)$ has a unique expression
\[
 \phi=\tau_\gamma\boldsymbol\phi_\gamma,\qquad
 \boldsymbol\phi_\gamma=(\phi_{\gamma,1},\dots,\phi_{\gamma,d})^t\in\mcS(\Omega_\gamma)^d.
\]
For $\mathbf u_\gamma=(u_{\gamma,1},\dots,u_{\gamma,d})^t\in\mcS'(\Omega_\gamma)^d$, the notation $u=\tau_\gamma^\vee\mathbf u_\gamma$ denotes the functional
\[
 \langle u,\phi\rangle=\langle\mathbf u_\gamma,\boldsymbol\phi_\gamma\rangle
 :=\sum_{i=1}^d\langle u_{\gamma,i},\phi_{\gamma,i}\rangle.
\]
Every $u\in\mcS'(\Omega_\gamma,\mcE_\lambda)$ is represented uniquely in this way. We call $\tau_\gamma$ the \emph{canonical local frame}, or \emph{orbit frame}, of $\mcE_\lambda|_{\Omega_\gamma}$.

The $F_w$-coordinate of $\Omega_\gamma\cong B_\gamma\times F_w$ depends on $\gamma$. At a point of $\Omega_\gamma\cap\Omega_{\gamma'}$, the coordinates $(b,f)$ and $(b,f')$ satisfy $f'=\Ad_{s_{\gamma'}(b)^{-1}s_\gamma(b)}f$. Define the transition between the orbit lifts and the corresponding tempered bundle transition matrix by
\[
 T_\gamma^{\gamma'}(b,f)=g_{\gamma'}(b,f')^{-1}g_\gamma(b,f)\in P_0,\qquad
 A_\gamma^{\gamma'}=\rho_{\sigma,\lambda}(T_\gamma^{\gamma'}).
\]
Thus, on the overlap,
\begin{equation}\label{eq:tau-transition}
 \tau_\gamma=\tau_{\gamma'}A_\gamma^{\gamma'},\qquad
 \tau_\gamma^\vee=\tau_{\gamma'}^\vee(A_\gamma^{\gamma'})^{-t}.
\end{equation}
The coefficient columns of the same section or distribution consequently satisfy
\[
 \boldsymbol\phi_{\gamma'}=A_\gamma^{\gamma'}\boldsymbol\phi_\gamma,\qquad
 \mathbf u_{\gamma'}=(A_\gamma^{\gamma'})^{-t}\mathbf u_\gamma.
\]
These equalities are at the same point of the overlap in the respective coordinates.

\begin{lemma}[$N$-invariance of the orbit frames]\label{lem:tau-invariant}
The frame $\tau_\gamma$ is $N$-invariant. Under $\mcS(\Omega_\gamma,\mcE_\lambda)\cong\mcS(B_\gamma\times F_w)\otimes V_{\sigma,\lambda}$, the group $N$ acts by the twisted geometric action $n\star(b,f)=(b,\Ad_{s_\gamma(b)^{-1}}(n)f)$ on the first factor and trivially on the second. Dually, $\xi\in\fn$ acts on $\mcS'(\Omega_\gamma,\mcE_\lambda)\cong\mcS'(B_\gamma\times F_w)\otimes V_\sigma^*$ through
\[
 (\xi^\star\phi)(b,y)=\frac{\dd}{\dd t}\Big|_{t=0}\phi\bigl(b,\exp(-t\Ad(s_\gamma(b)^{-1})\xi)\,y\bigr)\qquad(\phi\in\mcS(B_\gamma\times F_w))
\]
and trivially on $V_\sigma^*$. (The trivialization itself does not depend on the choice of basis.)
\end{lemma}

\begin{proof}
For $n\in N$ and $x=s_\gamma(b)u\dot wP_0$ we have $n^{-1}x=s_\gamma(b)\Ad_{s_\gamma(b)^{-1}}(n^{-1})u\dot wP_0$, where $\Ad_{s_\gamma(b)^{-1}}(n^{-1})u\in N$ because $s_\gamma(b)\in L$ normalizes $N$. Hence
\[
 (n\cdot\tau_\gamma)(x)=[ns_\gamma(b)\Ad_{s_\gamma(b)^{-1}}(n^{-1})u\dot w,\mathbf e]
 =[s_\gamma(b)u\dot w,\mathbf e]=\tau_\gamma(x).
\]
So in these frames $N$ acts only on the coefficient columns, componentwise through its action on $\Omega_\gamma$, which in the coordinates $(b,f)$ is $n\star(b,f)$.
\end{proof}

\subsection{Polynomial densities on the fibre}\label{ssec:poly-densities}

Let $\Omega$ be a nonzero $\U$-invariant algebraic top form on $\F_w=\U/\U_w$ and $\mu=|\Omega|$ be the associated $N$-invariant analytic density on $F_w$, and put
\[
 \mcP(F_w)=\mathcal O(\F_w)_\C\,\mu,\qquad\mcP_k(F_w)=\{\omega\in\mcP(F_w):\IN^k\omega=0\}.
\]

\begin{proposition}[regularity along the fibre]\label{prop:poly-densities}
For every $k\ge1$, $\mcS'(F_w)^{\IN^k}=\mcP_k(F_w)$ is finite dimensional,
and $\mcS'(F_w)^{[\fn]}=\mcP(F_w)$. In the product chart $\Omega_\gamma\cong B_\gamma\times F_w$, with the twisted action of Lemma~\ref{lem:tau-invariant},
\[
 \mcS'(B_\gamma\times F_w)^{\IN^k}=\mcS'(B_\gamma)\otimes\mcP_k(F_w).
\]
Thus in the orbit frames, $\mcS'(\Omega_\gamma,\mcE_\lambda)^{\IN^k}=\mcS'(B_\gamma)\otimes\mcP_k(F_w)\otimes V_\sigma^*$.
\end{proposition}

\begin{proof}
Set $H=N\cap\dot wN_0^-\dot w^{-1}$, which is identified with $F_w$ by $c_w$ from Remark~\ref{rem:cw}. By \cite[Proposition~A.4]{Abe12}, $\IN^kT=0$ implies $T=f\mu$, where $f$ is polynomial in exponential coordinates of $H$, of degree bounded in terms of $k$. These coordinates and $c_w$ are algebraic, so $f\in\mathcal O(\F_w)_\C$.
$\mcP(F_w)$ is a locally finite rational $\U$-module, so every element is $\IN$-finite by Lemma~\ref{lem:commutator}(b). Hence $\mcS'(F_w)^{[\fn]}=\mcP(F_w)$.

The product version is immediate since the twisted action differentiates only in the fibre variable, and the twisted $\IN^k$-annihilation is equivalent to ordinary fibrewise $\IN^k$-annihilation.
\end{proof}

\section{An equivariant resolution of the fibre closure}\label{sec:resolution}

This section prepares the extension of polynomial densities across the boundary of the closure of $F_w$. It produces three things: the proper $Q$-equivariant map $p_Z:\tilde Z_w\to X$ with image $\overline{O_w}$, which is an isomorphism over $O_w$ (Lemma~\ref{lem:closure}); frames $\epsilon_{\gamma,v}$ of $p_Z^*\mcE_\lambda$ that extend across the boundary (Lemma~\ref{lem:trivial-pZE}); and the comparison $\epsilon_{\gamma,v}=\tau_\gamma A_v$ with the orbit frames, together with the formula for $A_v$ \eqref{eq:positive-Delta-order} (Lemma~\ref{lem:comparison}), which is the input for the pole analysis of Section~\ref{sec:regularized}.

\subsection{Resolution}\label{ssec:resolution}

We use \emph{simple normal crossings} (SNC) in the sense of \cite[Definition~24]{Kol07}. A reduced divisor $D$ on a smooth variety $Z$ over a field $k$ is SNC if its irreducible components are smooth over $k$ and, at each point $x\in Z$, the local equations of the components through $x$ form part of a regular system of parameters of $\mathcal O_{Z,x}$. In particular, at a real point of a smooth $\R$-variety these equations can be completed to an \'etale coordinate map defined over $\R$.

Recall $\F_w=\U\dot w\pP_0/\pP_0,\ \rR_w=\lL\cap\dot w\pP_0\dot w^{-1}$. By \eqref{eq:N-orbit-real-points}, $\F_w(\R)=N\dot wP_0/P_0=F_w\cong N/N_w$. Let $\Y_w\subset\X$ be the scheme-theoretic closure of $\F_w$ in $\X$; it is a projective variety over $\R$. Since $\F_w$ is stable under $\rR_w\ltimes\U$, so is $\Y_w$. The following theorem is obtained by applying functorial resolution to $\Y_w$ and then principalizing the boundary ideal on the resolved variety \cite[Theorems~35--36]{Kol07}. The $\rR_w\ltimes\U$-action lifts by smooth functoriality; see \cite[Proposition~9.1]{Kol07}.

\begin{theorem}[resolution of singularities]\label{thm:resolution}
There is an $\rR_w\ltimes\U$-equivariant resolution $\pi_Y:\widetilde{\Y_w}\to\Y_w$ such that
\begin{enumerate}[label=\textup{(\roman*)}]
\item $\widetilde{\Y_w}$ is a smooth projective variety over $\R$ with an action of $\rR_w\ltimes\U$;
\item $\pi_Y$ is $\rR_w\ltimes\U$-equivariant, projective and birational, and an isomorphism over $\F_w$;
\item $\D_w:=(\widetilde{\Y_w}\setminus\F_w)_{\rm red}$ is an $\rR_w\ltimes\U$-stable SNC divisor over $\R$.
\end{enumerate}
\end{theorem}

Since $\oO_w\cong\lL\times^{\rR_w}\F_w$, we put
\[
 \zZ_w=\lL\times^{\rR_w}\Y_w,\qquad\tilde\zZ_w=\lL\times^{\rR_w}\widetilde{\Y_w},\qquad\boldsymbol\Delta_w=(\tilde\zZ_w\setminus\oO_w)_{\rm red}\cong\lL\times^{\rR_w}\D_w.
\]
These are projective varieties over $\R$, and $\tilde\zZ_w$ is smooth; $\boldsymbol\Delta_w$ is an SNC divisor over $\R$, as seen in the local charts of the associated bundle. Since $\rR_w\ltimes\U$ is connected, it preserves each component of $\D_w$; similarly, $\Q$ preserves each component of $\boldsymbol\Delta_w$.
We have the diagram
\begin{equation}\label{eq:resolution-diagram}
\begin{tikzcd}
 \lL\times^{\rR_w}\F_w\arrow[r,hook]\arrow[d,"\simeq"'] & \tilde\zZ_w=\lL\times^{\rR_w}\widetilde{\Y_w}\arrow[d,"\pi_Z"']\arrow[rdd,bend left,"p_Z"] & \\
 \lL\times^{\rR_w}\F_w\arrow[r,hook]\arrow[d,"\simeq"'] & \zZ_w=\lL\times^{\rR_w}\Y_w\arrow[d,"m_Z"'] & \\
 \oO_w\arrow[r,hook] & \overline{\oO_w}\arrow[r,hook] & \X
\end{tikzcd}
\end{equation}
where $\pi_Z$ is induced by $\pi_Y$, the action map $\lL\times\Y_w\to\overline{\oO_w}$ factors through $m_Z$, which is proper, surjective and birational, and $p_Z=m_Z\circ\pi_Z$ is proper. The maps $\pi_Z$, $m_Z$ and $p_Z$ are $\Q$-equivariant.

\begin{lemma}\label{lem:real-points}
On real points,
\[
\begin{aligned}
 (\lL\times^{\rR_w}\F_w)(\R)&\cong\lL(\R)\times^{\rR_w(\R)}\F_w(\R)\cong L\times^{R_w}F_w\cong O_w,\\
 (\lL\times^{\rR_w}\Y_w)(\R)&\cong\lL(\R)\times^{\rR_w(\R)}\Y_w(\R)\cong L\times^{R_w}Y_w=:Z_w,\\
 (\lL\times^{\rR_w}\widetilde{\Y_w})(\R)&\cong\lL(\R)\times^{\rR_w(\R)}\widetilde{\Y_w}(\R)\cong L\times^{R_w}\widetilde{Y_w}=:\tilde Z_w,
\end{aligned}
\]
where $Y_w:=\Y_w(\R)\subset X$ is a closed semialgebraic subset and $\widetilde{Y_w}:=\widetilde{\Y_w}(\R)$ is a connected projective Nash manifold.
\end{lemma}

\begin{proof}
For $\mathbf S\in\{\F_w,\Y_w,\widetilde{\Y_w}\}$, let
$p:\lL\times^{\rR_w}\mathbf S\longrightarrow\lL/\rR_w$
be the projection. Take $x\in(\lL\times^{\rR_w}\mathbf S)(\R)$ and put $b=p(x)$. By Proposition~\ref{prop:flag-geometry} applied to $\lL$, choose $l\in L$ representing $b$. Since $l$ is a real point, the map
$\mathbf S\xrightarrow{\sim}p^{-1}(b),\ s\longmapsto[l,s]$
is an isomorphism over $\R$. Thus $x=[l,s]$ for some $s\in\mathbf S(\R)$, proving surjectivity of the natural map
$L\times^{R_w}\mathbf S(\R)\longrightarrow(\lL\times^{\rR_w}\mathbf S)(\R)$.

If $[l,s]=[l',s']$ with $l,l'\in L$ and $s,s'\in\mathbf S(\R)$, then
$r=l^{-1}l'\in L\cap\rR_w(\R)=R_w,\ s'=r^{-1}s$,
so the map is also injective.
By Theorem~\ref{thm:resolution}, $\widetilde{\Y_w}(\R)$ is a projective Nash manifold. It is connected because $F_w$ is dense and connected.
\end{proof}
 
\begin{lemma}\label{lem:closure}
We have $\Y_w(\R)=\overline{F_w}$. Consequently,
\begin{equation}\label{eq:pZ-open-preimage}
 p_Z(\tilde Z_w)=\overline{O_w},\qquad p_Z^{-1}(O_w)=L\times^{R_w}F_w,
\end{equation}
and $p_Z$ is an isomorphism over $O_w$.
\end{lemma}

\begin{proof}
Since $w\in{}_IW$, we have $\F_w=\pP_0\dot w\pP_0/\pP_0$. The relative Bruhat closure theorem \cite[Theorem~3.13 and Corollary~3.15]{BT72}, passed to the flag variety, gives $Y_w=\overline{F_w}=\bigsqcup_{v\le w}P_0\dot vP_0/P_0$. 
As $w$ is minimal in $W_Iw$, this implies $Y_w\cap O_w=F_w$. 
Since $\pi_Y$ is an isomorphism over $F_w$, we obtain $p_Z^{-1}(O_w)=L\times^{R_w}F_w$, on which $p_Z$ is the isomorphism $[l,f]\mapsto lf$. Finally, properness of $p_Z$ and density of $L\times^{R_w}F_w$ in $\tilde Z_w$ give $p_Z(\tilde Z_w)=\overline{O_w}$.
\end{proof}

Taking real points in \eqref{eq:resolution-diagram} and using the preceding lemmas, we obtain
\[
\begin{tikzcd}
 L\times^{R_w}F_w\arrow[r,hook]\arrow[d,"\simeq"'] & \tilde Z_w=L\times^{R_w}\widetilde{Y_w}\arrow[d,"\pi_Z"']\arrow[rdd,bend left,"p_Z"] & \\
 L\times^{R_w}F_w\arrow[r,hook]\arrow[d,"\simeq"'] & Z_w=L\times^{R_w}Y_w\arrow[d,"m_Z"'] & \\
 O_w\arrow[r,hook] & \overline{O_w}\arrow[r,hook] & X
\end{tikzcd}
\]
and $p_Z$ is a $Q$-equivariant proper Nash map between the projective Nash manifolds $\tilde Z_w$ and $X$.

\subsection{Adapted charts and frames}\label{ssec:charts-frames}

In this section, we trivialize the pullback bundle on $\tilde Z_w$, and define the local chart in which we define local extensions and meromorphic continuation.

Put $D_w=\widetilde{Y_w}\setminus F_w\cong\D_w(\R)$ and $\Delta_w=\tilde Z_w\setminus O_w\cong\boldsymbol\Delta_w(\R)\cong L\times^{R_w}D_w$. By the SNC convention of Section~\ref{ssec:resolution}, every real point $x\in\widetilde{Y_w}$ has an affine open neighbourhood $\V_x$ with an \'etale morphism
\[
 \boldsymbol\phi_x=(\phi_1,\dots,\phi_n):\V_x\to\mathbb A^n,\qquad\boldsymbol\phi_x^{-1}(\{t_i=0\})=\D_{w,i}\cap\V_x\quad(1\le i\le r),
\]
where $t_1,\dots,t_n$ are the standard coordinates on $\mathbb A^n$ and $\D_{w,1},\dots,\D_{w,r}$ are the irreducible components of $\D_w$ through $x$. Since $\boldsymbol\phi_x$ is \'etale, $\Omega^1_{\V_x/\R}$ is free on $\dd\phi_1,\dots,\dd\phi_n$. By the Nash inverse function theorem, after shrinking, there is a Nash open neighbourhood $V_x\subset\V_x(\R)$ of $x$ on which $\boldsymbol\phi_x$ restricts to a Nash open immersion $V_x\to\R^n$ with $\boldsymbol\phi_x(x)=0$ and
\[
 \boldsymbol\phi_x(D_{w,i}\cap V_x)=\{t\in\boldsymbol\phi_x(V_x):t_i=0\}\qquad(D_{w,i}=\D_{w,i}(\R),\ 1\le i\le r).
\]
Such Nash open sets $V_x$ form a basis of neighbourhoods of $x$ in the analytic topology.

To trivialize $p_Z^*\mcE_\lambda$, let $p_Y:\widetilde{\Y_w}\to\X$ be $\pi_Y$ composed with the closed immersion $\Y_w\subset\X$, and put
\[
 \cC_v=p_Y^{-1}(\dot v\U_0^-\pP_0/\pP_0),\qquad C_v=\cC_v(\R)=p_Y^{-1}(\dot vN_0^-P_0/P_0)\qquad(v\in W).
\]
By Proposition~\ref{prop:flag-geometry}, the sets $C_v$ cover $\widetilde{Y_w}$.
The following is immediate:

\begin{lemma}\label{lem:trivial-Bgamma}
The section $s_\gamma:\B_\gamma\to\lL$ gives algebraic trivializations of the associated bundles $\tilde\zZ_w\to\lL/\rR_w$ and $\zZ_w\to\lL/\rR_w$ over $\B_\gamma$. Namely, the map
\[
 \B_\gamma\times\widetilde{\Y_w}\to\tilde\zZ_w|_{\B_\gamma},\qquad(b,y)\mapsto[s_\gamma(b),y],
\]
is an isomorphism of varieties, with inverse $[l,y]\mapsto\bigl(l\rR_w,(s_\gamma(l\rR_w)^{-1}l)\cdot y\bigr)$, and similarly for $\zZ_w|_{\B_\gamma}$.
\end{lemma}

Let $\oOm_{\gamma,v}=\B_\gamma\times\cC_v$, regarded as a Zariski open subset of $\tilde\zZ_w|_{\B_\gamma}$ via Lemma~\ref{lem:trivial-Bgamma}, and let $\Omega_{\gamma,v}=\oOm_{\gamma,v}(\R)$, a Nash open subset of $\tilde Z_w|_{B_\gamma}$.

\begin{lemma}[trivialization of $p_Z^*\mcE_\lambda$]\label{lem:trivial-pZE}
Let $q_v:\cC_v\to\dot v\U_0^-$ be $p_Y|_{\cC_v}$ followed by $\dot v\U_0^-\pP_0/\pP_0\cong\dot v\U_0^-$, so that $q_v(y)\pP_0=p_Y(y)$. Define the \emph{resolution lift}
\[
 q_{\gamma,v}:\B_\gamma\times\cC_v\to\G,\qquad(b,y)\mapsto s_\gamma(b)\,q_v(y).
\]
It satisfies $q_{\gamma,v}(b,y)\pP_0=p_Z(b,y)$. Its restriction $q_{\gamma,v}:B_\gamma\times C_v\to G$ defines an isomorphism of tempered bundles
\[
 T_{\gamma,v}:\Omega_{\gamma,v}\times V_{\sigma,\lambda}\to(p_Z^*\mcE_\lambda)|_{\Omega_{\gamma,v}},\qquad\bigl((b,y),e\bigr)\mapsto\bigl((b,y),[q_{\gamma,v}(b,y),e]\bigr),
\]
with inverse $(z,[g,e])\mapsto(z,\rho_{\sigma,\lambda}(q_{\gamma,v}(z)^{-1}g)e)$.
\end{lemma}

\begin{proof}
The representatives $\dot\gamma$ and $\dot v$ were chosen in $G$, and all real points of unipotent groups lie in $G$, hence $q_{\gamma,v}(B_\gamma\times C_v)\subset G$. The lifts are Nash, and the trivialization is tempered by \cite[Section~6.1]{CS21}.
\end{proof}

With the basis of Section~\ref{ssec:orbit-frames}, the resolution lift defines the row frame
\[
 \epsilon_{\gamma,v}=(\epsilon_{\gamma,v,1},\dots,\epsilon_{\gamma,v,d}),\qquad
 \epsilon_{\gamma,v,j}(b,y)=\bigl((b,y),[q_{\gamma,v}(b,y),e_j]\bigr).
\]
Every $\Phi\in\mcS(\Omega_{\gamma,v},p_Z^*\mcE_\lambda)$ has a unique expression $\Phi=\epsilon_{\gamma,v}\boldsymbol\Phi$, where $\boldsymbol\Phi=(\Phi_1,\dots,\Phi_d)^t\in\mcS(\Omega_{\gamma,v})^d$. We call $\epsilon_{\gamma,v}$ the \emph{canonical local frame}, or \emph{resolution frame}, on $\Omega_{\gamma,v}$ and write $\epsilon_{\gamma,v}^\vee$ for its dual row frame. At a point of $\Omega_{\gamma,v}\cap\Omega_{\gamma',v'}$, put $y'=s_{\gamma'}(b)^{-1}s_\gamma(b)y$ and define the transition between the resolution lifts
\[
 R^{\gamma',v'}_{\gamma,v}(b,y)=q_{\gamma',v'}(b,y')^{-1}q_{\gamma,v}(b,y)\in P_0.
\]
Its image under $\rho_{\sigma,\lambda}$ is the bundle transition matrix. Thus, writing $\boldsymbol\Phi'$ for the coefficient column in the frame $\epsilon_{\gamma',v'}$, we have
\begin{equation}\label{eq:eps-transition}
 \epsilon_{\gamma,v}=\epsilon_{\gamma',v'}\rho_{\sigma,\lambda}(R^{\gamma',v'}_{\gamma,v}),\qquad
 \boldsymbol\Phi'=\rho_{\sigma,\lambda}(R^{\gamma',v'}_{\gamma,v})\boldsymbol\Phi.
\end{equation}

Since the $C_v$ cover $\widetilde{Y_w}$, every real point $x\in\widetilde{Y_w}$ has an affine open neighbourhood $\V_x$ as above contained in some $\cC_v$, together with a Nash open $V_x\subset\V_x(\R)$ as above. We call the quadruple $(\cC_v,\V_x,V_x,\boldsymbol\phi_x)$ an \emph{adapted chart} at $x$, with respect to the SNC divisor $\D_w$.

\subsection{The comparison factor}\label{ssec:comparison}

We emphasize that the polynomial density description along the fibres uses the orbit frames $\tau_\gamma$, which need not extend across the boundary. We therefore compare them with the frames $\epsilon_{\gamma,v}$ defined across the boundary of $\tilde Z_w$. On their common domain, the two frames differ by a tempered transition depending on $\lambda$. This matrix enters the local regularized pairings and requires careful consideration in meromorphic continuation, so we describe it explicitly below.

\begin{definition}[comparison factor]\label{def:comparison}
Let $\cC_v^\circ=\cC_v\setminus\D_w\subset\F_w$ and $C_v^\circ=\cC_v^\circ(\R)=C_v\setminus D_w\subset F_w$. On $\B_\gamma\times\cC_v^\circ$, the orbit lift $g_\gamma$ and the resolution lift $q_{\gamma,v}$ have the same projection to $\X$, namely $p_Z$. Define
\[
 \lLd^\gamma_v:\B_\gamma\times\cC_v^\circ\to\pP_0,\qquad
 \lLd^\gamma_v=g_\gamma^{-1}q_{\gamma,v},
\]
let $\overline{\lLd^\gamma_v}=\pr_{\lL_0}\circ\lLd^\gamma_v$ be its Levi component, and let $\Lambda^\gamma_v:B_\gamma\times C_v^\circ\to P_0$ and $\overline{\Lambda^\gamma_v}:B_\gamma\times C_v^\circ\to L_0$ be the maps on real points. Put
\[
 A^\gamma_v=\rho_{\sigma,\lambda}(\Lambda^\gamma_v)=\rho_{\sigma,\lambda}(\overline{\Lambda^\gamma_v}),
\]
an invertible $d\times d$ matrix of tempered functions on $B_\gamma\times C_v^\circ$, with respect to the basis $e_1,\dots,e_d$.
\end{definition}

Expanding the two lifts shows that their common base factor $s_\gamma(b)$ cancels:
\begin{equation}\label{eq:base-lift-cancels}
\begin{aligned}
 \Lambda^\gamma_v(b,y)
 =\bigl(s_\gamma(b)c_w(y)\dot w\bigr)^{-1}s_\gamma(b)q_v(y)
 =(c_w(y)\dot w)^{-1}q_v(y).
\end{aligned}
\end{equation}
Thus $\Lambda^\gamma_v$, its Levi component and $A^\gamma_v$ do not depend on $b$ or $\gamma$ in these product charts, and we write $\Lambda_v$, $\overline\Lambda_v$ and $A_v$. On $B_\gamma\times C_v^\circ$, identifying the orbit frame $\tau_\gamma$ with its pullback by $p_Z$ we obtain
\begin{equation}\label{eq:eps-tau}
 \epsilon_{\gamma,v}=\tau_\gamma A_v,\qquad
 \epsilon_{\gamma,v}^\vee=\tau_\gamma^\vee A_v^{-t}.
\end{equation}

The following lemma separates the dependence of $A_v$ on $\lambda$ into powers of real characters, and shows that one of these characters vanishes on the whole boundary.

\begin{lemma}[the comparison factor]\label{lem:comparison}
Put $\Delta_0(p)=\det(\Ad(p)|_{\Lie\U_0})$ for $p\in\pP_0$. Choose real algebraic characters $\chi_1,\dots,\chi_{n_0}$ of $\lL_0$, where $n_0=\dim\fa_0$, with $\chi_1=\Delta_0|_{\lL_0}$, whose differentials on $\fa_0$ form a real basis of $\fa_0^*$, and write $\lambda=\sum_k\lambda_k\,\dd\chi_k|_{\fa_0}, \ \lambda_k\in \C$. Then
\[
 A_v=\left(\prod_{k=1}^{n_0}|\chi_k\circ\overline\Lambda_v|^{\lambda_k}\right)m_v,
\]
where $m_v$ is a bounded Nash matrix, independent of the base variable and of $\lambda$. Moreover, every boundary component $D$ meeting $C_v$ satisfies
\begin{equation}\label{eq:positive-Delta-order}
 \ord_D(\Delta_0\circ\overline\Lambda_v)>0.
\end{equation}
\end{lemma}

\begin{proof}
Write $\overline\Lambda_v(y)=m(y)a(y)$ with $m(y)\in M_0$ and $a(y)\in A_0$. Compactness of $M_0$ gives $|\chi_k(m(y))|=1$, while $\log|\chi_k(a(y))|=\dd\chi_k(\log a(y))$. Hence
\[
 \rho_{\sigma,\lambda}(\overline\Lambda_v(y))=e^{\lambda(\log a(y))}\sigma(m(y))=\prod_k|\chi_k(\overline\Lambda_v(y))|^{\lambda_k}\,\sigma(m(y)),
\]
which is the formula with $m_v=\sigma(m(y))$. The Langlands projections are Nash and $\sigma$ is Nash, so $m_v$ is a Nash matrix; compactness of $M_0$ gives boundedness.

Set $\mcL=p_Y^*(\det T_\X)$ line bundle on $\widetilde{\Y_w}$. Choose a basis $E_1,\dots,E_m$ of $\Lie\U_0^-$, $m=\dim\X$, and let $\bar E_i$ be its image in $T_o\X=\fg/\fp_0$, with $o=e\pP_0$. We have
$\det T_\X\simeq\G\times^{\pP_0}\mathbb A^1_{\Delta_0^{-1}}$, and fix $\bar E_1\wedge\cdots\wedge\bar E_m\mapsto [e,1]$ at base point $o$.
Since $q_v(y)\pP_0=p_Y(y)$, the algebraic trivialization
\[
 \cC_v\times\mathbb A^1\xrightarrow{\sim}\mcL|_{\cC_v},
 \qquad (y,t)\longmapsto(y,[q_v(y),t])
\]
defines the nowhere-vanishing frame $e_v(y)=(y,[q_v(y),1])=(y,q_v(y)\cdot[e,1])$ on all of $\cC_v$.

Put $\xi_i=\Ad(\dot w)E_i$ and $s_w=\xi_{1,\X}\wedge\cdots\wedge\xi_{m,\X}$, where $\xi_\X(x)=\frac{\dd}{\dd t}|_{t=0}\exp(t\xi)x$. Write $p_Y^*s_w=f_ve_v$, with $f_v\in\mathcal O(\cC_v)$. Translating tangent vectors at $q_v(y)\pP_0$ back to $o$ gives
\[
 f_v(y)\,\bar E_1\wedge\cdots\wedge\bar E_m
 =\bigwedge_{i=1}^m\overline{\Ad(q_v(y)^{-1})\xi_i}.
\]
On $\cC_v^\circ$, let $\lLd_v=(c_w\dot w)^{-1}q_v$ be the algebraic comparison map. With $u=c_w(y)$, we have
\[
 q_v(y)=u\dot w\,\lLd_v(y),\qquad
 \Ad(q_v(y)^{-1})\xi_i
 =\Ad(\lLd_v(y)^{-1})\Ad(\dot w^{-1}u^{-1}\dot w)E_i.
\]
The second factor acts on $\Lie\U_0^-$ with determinant one, since $\dot w^{-1}u^{-1}\dot w\in\U_0^-$. The first factor induces an operator on $\fg/\fp_0$ with determinant $\Delta_0^{-1}(\lLd_v(y)^{-1})=\Delta_0(\lLd_v(y))$. Thus
\[
 f_v|_{\cC_v^\circ}=\Delta_0\circ\lLd_v
 =\Delta_0\circ\overline{\lLd_v}.
\]

Finally, $s_w$ is nonzero precisely on the unique open $\dot w\U_0^-\dot w^{-1}$-orbit $\dot w\U_0^-\pP_0/\pP_0$, since its factors span the orbit tangent space. As $\Y_w\cap(\dot w\U_0^-\pP_0/\pP_0)=\F_w$, the section $p_Y^*s_w$ vanishes precisely on $\D_w$, so $f_v$ has positive order along every boundary component meeting $\cC_v$.
\end{proof}

\section{Meromorphic continuation of Nash integrals}\label{sec:regularized}

This section gives a continuation theorem for Mellin integrals with bounded Nash coefficients. We then apply it to the local regularized pairings and obtain a common set of good parameters for all polynomial densities. The proof of the continuation theorem is deferred to Section~\ref{ssec:monomial-proof}.

The continuation is obtained by subtracting boundary Taylor polynomials and integrating their monomials explicitly, following the resolution approach of Atiyah~\cite{Ati70} and the rectilinearization methods of Bierstone--Milman~\cite{BM88} and Parusi\'nski~\cite{Par94}. A bounded Nash coefficient need not extend smoothly to the boundary. We therefore resolve the Zariski closure of its graph and principalize the boundary data, obtaining smooth amplitudes and monomial powers on closed coordinate orthants.

\subsection{Mellin continuation with Nash coefficients}

Put $C=(0,1)^n$, fix $0\le r\le n$, and let $m:C\to\C$ be a bounded Nash function, meaning that its real and imaginary parts are Nash. For a Nash manifold $B$ and $s=(s_1,\ldots,s_r)\in\C^r$ with $\Re s_i>-1$, define
\begin{equation}\label{eq:Nash-integral}
 (T(s)\Phi)(b)=\int_C m(x)\Phi(b,x)\prod_{i=1}^r x_i^{s_i}\,\dd x,
 \qquad \Phi\in\mcS(B\times\R^n).
\end{equation}
Here $x_i^{s_i}=\exp(s_i\log x_i)$. Boundedness of $m$ gives absolute convergence and a holomorphic family of continuous maps to $\mcS(B)$.

\begin{theorem}[Mellin continuation with Nash coefficients]\label{thm:Nash-continuation}
The operators $T(s):\mcS(B\times\R^n)\to\mcS(B)$ extend meromorphically to $\C^r$. There are finitely many vectors $a_\nu\in\mathbb Z_{\ge0}^r\setminus\{0\}$, depending only on $m$ and $r$, such that all poles are contained in the locally finite collection of hyperplanes
\begin{equation}\label{eq:Nash-poles}
 \langle a_\nu,s\rangle+j+1=0,\qquad j\in\mathbb Z_{\ge0}.
\end{equation}
For every compact $K\subset\C^r$, there is a polynomial $H_K$, a finite product of the affine-linear functions $\langle a_\nu,s\rangle+j+1$, for which the maps
\[
 s\longmapsto H_K(s)T(s)\Phi,\qquad \Phi\in\mcS(B\times\R^n),
\]
are holomorphic with values in $\mcS(B)$ on a common open neighbourhood of $K$. For every continuous seminorm $p$ on $\mcS(B)$ there are a continuous seminorm $q$ on $\mcS(B\times\R^n)$ and a constant $C_K$ such that
\begin{equation}\label{eq:Mellin-estimate}
 \sup_{s\in K}p\bigl(H_K(s)T(s)\Phi\bigr)\le C_Kq(\Phi)
 \qquad\bigl(\Phi\in\mcS(B\times\R^n)\bigr).
\end{equation}
\end{theorem}

The same continuation holds for holomorphic $\mcS(B\times\R^n)$-valued amplitudes $\Phi(s)$, with same candidate poles. By \eqref{eq:Mellin-estimate}, for every fixed $u\in\mcS'(B)$, the functionals $\Phi\mapsto\langle u,T(s)\Phi\rangle$ form a meromorphic family in $\mcS'(B\times\R^n)$, equipped with its strong dual topology.

Let $V:\C^p\to\C^r$ be a complex-linear map defined over $\R$, and suppose there is $v\in\R^p$ such that every coordinate of $Vv$ is positive. Put
\[
 \ell_\nu(\lambda)=\langle a_\nu,V\lambda\rangle,
 \qquad
 \mathcal H_V=\bigcup_\nu\ \bigcup_{N\in\mathbb Z}
 \{\lambda\in\C^p:\ell_\nu(\lambda)=N\}.
\]
When there are no candidate poles, this union is empty.
Note that $\ell_\nu(v)=\langle a_\nu,Vv\rangle>0$, so each hyperplane is proper and transverse to $v$.
The following is an easy consequence of Theorem~\ref{thm:Nash-continuation}:

\begin{lemma}\label{lem:integer-shifts}
The set $\mathcal H_V$ is a locally finite union of proper affine hyperplanes, all transverse to $v$. For every $d\in\mathbb Z^r$ and every $\lambda\notin\mathcal H_V$, the meromorphic family
\[
 T(d+z+V\lambda),\qquad z\in\C^r,
\]
is holomorphic at $z=0$, and defines a continuous map $\mcS(B\times\R^n)\to\mcS(B)$. The same assertion holds for holomorphic family $\Phi(z)$.
\end{lemma}

\subsection{Application to the local regularized pairings}\label{ssec:local-pairing}

Cover $\widetilde{Y_w}$ by finitely many adapted charts $(\cC_v,\V_x,V_x,\boldsymbol\phi_x)$, with $\boldsymbol\phi_x=(\phi_1,\ldots,\phi_n)$ and $\D_w\cap\V_x=\{\phi_1\cdots\phi_r=0\}$. Take $V_x=(-\delta,\delta)^n$ in these coordinates, with closure contained in a larger adapted coordinate neighbourhood. In the canonical local frame $\epsilon_{\gamma,v}$ we identify test sections with coefficient columns $\boldsymbol\Phi\in\mcS(B_\gamma\times V_x)^d$.

The input polynomial densities are expressed in the orbit frame $\tau_\gamma$, while the test sections use the resolution frame $\epsilon_{\gamma,v}$, which extends across the boundary. Thus the comparison matrix $A_v$ in \eqref{eq:eps-tau} enters the fibre integral. For a polynomial density $\omega\in\mcP(F_w)$, initially put
\begin{equation}\label{eq:local-fibre-integral}
 J_\omega(z,\lambda;\boldsymbol\Phi)(b)
 =\int_{V_x^\circ}A_v(y)\boldsymbol\Phi(b,y)\,\omega(y)
        \prod_{l=1}^r|\phi_l(y)|^{z_l},
 \qquad V_x^\circ=V_x\cap F_w,
\end{equation}
for $\Re z_l$ sufficiently large. We apply the general theorem to this integral before taking its value at $z=0$.

\emph{The exponents and the smooth amplitude.}
For $\omega=f\mu\ne0$, write on $\V_x$
\begin{equation}\label{eq:omega-decomp}
 \Omega=h\prod_{l=1}^r\phi_l^{e_l}\,\dd\phi_1\wedge\cdots\wedge\dd\phi_n,
 \qquad f=f_0\prod_{l=1}^r\phi_l^{b_l},
 \qquad d_{\omega,l}=e_l+b_l\in\mathbb Z.
\end{equation}
Here $e_l$ and $b_l$ are the divisor orders, positive for zeros and negative for poles. After extracting these powers, normality gives a regular function $f_0$ and a regular unit $h$ on $\V_x$; $f_0$ may vanish and may be complex-valued. On an orthant with $\operatorname{sgn}\phi_l=\varepsilon_l$, the density is
\[
 \omega=f_\varepsilon\prod_{l=1}^r|\phi_l|^{d_{\omega,l}}\,|\dd\phi|,
 \qquad f_\varepsilon=f_0|h|\prod_{l=1}^r\varepsilon_l^{b_l}.
\]
Likewise,
\begin{equation}\label{eq:nu-kl}
 \chi_k\circ\overline\Lambda_v=g_k\prod_{l=1}^r\phi_l^{\nu_{kl}},
 \qquad \nu_{kl}\in\mathbb Z,
\end{equation}
where the $g_k$ are regular units on $\V_x$. Lemma~\ref{lem:comparison} gives $\nu_{1l}>0$ for every boundary index $l$.

Fix an orthant and reflect its coordinates to $C=(0,\delta)^n$. By Lemma~\ref{lem:comparison}, each scalar term of \eqref{eq:local-fibre-integral} has the form
\begin{equation}\label{eq:orthant-integral}
 \int_C m(y)\,\Psi_{b,\lambda}(y)\prod_{l=1}^r y_l^{s_l}\,\dd y,
 \qquad s_l=z_l+d_{\omega,l}+\sum_k\nu_{kl}\lambda_k,
\end{equation}
where $m$ is a scalar entry of the bounded Nash matrix $m_v$, and, before reflection,
\[
 \Psi_{b,\lambda}(y)=\Phi_j(b,y)\,f_\varepsilon(y)
                           \prod_k|g_k(y)|^{\lambda_k}.
\]
This multiplier is smooth on a neighbourhood of the closed coordinate cube, with all derivatives bounded locally uniformly in $\lambda$. Thus it gives a holomorphic continuous multiplier on $\mcS(B_\gamma\times V_x)$. First extend $\Psi$ by zero from the full cube $V_x$ to $\R^n$, using \cite[Theorem~5.4.1]{AG08}, and then reflect and restrict to $C$.

Under $y=\delta t$, formula~\eqref{eq:orthant-integral} becomes \eqref{eq:Nash-integral} with coefficient $m_\delta(t)=m(\delta t)$ and amplitude $\Psi_\delta(b,t;\lambda)=\Psi_{b,\lambda}(\delta t)$, multiplied by $\delta^{n+\sum_{l=1}^r s_l}$, a scalar factor nowhere-zero and entire in $s$. Applying Theorem~\ref{thm:Nash-continuation} gives the meromorphic continuation and Schwartz estimates. The candidate pole directions depend only on the chart, the orthant and the entries of $m_v$. All dependence on $\omega$ occurs in the integer vector $d_\omega$ and in the smooth amplitude, so the same candidate pole directions work for every polynomial density.

\emph{A common set of good parameters.}
In each chart set $(V\lambda)_l=\sum_k\nu_{kl}\lambda_k$. For each vector $a_\nu$ supplied by Theorem~\ref{thm:Nash-continuation}, put
\[
 \ell_\nu(\lambda)=\langle a_\nu,V\lambda\rangle
 =\sum_k\left(\sum_l a_{\nu,l}\nu_{kl}\right)\lambda_k.
\]
Let $\mathcal L_w$ be the finite list of these linear functions given by the fixed charts, orthants and matrix entries, and define
\[
 \mcH_w=\bigcup_{\ell\in\mathcal L_w}\ \bigcup_{N\in\mathbb Z}
                      \{\lambda\in\fa_{0,\C}^*:\ell(\lambda)=N\}.
\]

Note that for $\nu_0=\dd\Delta_0|_{\fa_0}$, we have $\lambda_1=1,\lambda_k=0$ for $k>1$, and every coordinate of $V\nu_0$ is $\nu_{1l}>0$. Apply Lemma~\ref{lem:integer-shifts} to the finitely many scalar integrals above gives the following:

\begin{theorem}[good parameters]\label{thm:good-parameters}
The set $\mcH_w$ is a locally finite union of proper affine hyperplanes. If
\begin{equation}\label{eq:good-parameter}
 \lambda\notin\mcH_w,
\end{equation}
then, for every polynomial density $\omega$, the family $J_\omega(z,\lambda;\boldsymbol\Phi)$ is holomorphic at $z=0$, and $\boldsymbol\Phi\mapsto J_\omega(0,\lambda;\boldsymbol\Phi)$ is a continuous map to $\mcS(B_\gamma)^d$. Such $\lambda$ are called \emph{good parameters}; the remaining parameters are called \emph{bad}.
\end{theorem}

\begin{remark}[a transverse direction]\label{rem:transverse}
The good parameters form an open set, and the same direction $\nu_0$ is transverse to every hyperplane of $\mcH_w$, since
$\ell_\nu(\nu_0)=\sum_l a_{\nu,l}\nu_{1l}>0$.
By local finiteness, for every $\lambda_0$ there is $\epsilon>0$ such that $\lambda_0+\tau\nu_0$ is good for $0<|\tau|<\epsilon$. Since $\nu_0$ is independent of $w$, the same holds simultaneously for the finitely many $Q$-orbits.
\end{remark}

\begin{definition}[local regularized pairing]\label{def:local-pairing}
Let $\lambda$ be good and let $u\in\mcS'(O_w,\mcE_\lambda)^{\IN^k}$. Localize using the fixed base partition chosen in Section~\ref{ssec:global-presentation}:
$u=\sum_\gamma u_\gamma$, where $u_\gamma=(\beta_\gamma\circ\pi_w)u$ is supported in $C_\gamma=\pi_w^{-1}(K_\gamma)$. Choose a basis $(\omega_p)$ of $\mcP_k(F_w)$. Proposition~\ref{prop:poly-densities} gives
\[
 u_\gamma=\tau_\gamma^\vee\sum_p(\mathbf u_p\otimes\omega_p),\qquad
 \mathbf u_p\in\mcS'(B_\gamma)^d,\qquad \supp\mathbf u_p\subset K_\gamma.
\]
For $\Phi=\epsilon_{\gamma,v}\boldsymbol\Phi\in\mcS(B_\gamma\times V_x,p_Z^*\mcE_\lambda)$, define
\begin{equation}\label{eq:local-pairing}
 \Bigl\langle\int^{\lambda,\gamma,v}u_\gamma,\Phi\Bigr\rangle
 =\sum_p\bigl\langle\mathbf u_p,J_{\omega_p}(0,\lambda;\boldsymbol\Phi)\bigr\rangle.
\end{equation}
Theorem~\ref{thm:good-parameters} makes this a (continuous) tempered distribution. For fixed coefficient distributions $\mathbf u_p$, the same formula is meromorphic in $\lambda$ by Theorem~\ref{thm:Nash-continuation}. Compatibility of these local pairings is proved in Section~\ref{sec:zero-order}.
\end{definition}

\subsection{Proof of Theorem~\ref{thm:Nash-continuation}}\label{ssec:monomial-proof}

\begin{lemma}[Mellin continuation]\label{lem:Mellin}
Let $\Omega\subset\C^p$ be connected and open, and let $A:\Omega\to C^\infty([0,1]^n)$ be holomorphic. Let $L_1,\ldots,L_n$ be linear functions on $\C^p$ such that
\[
 \Omega_+=\{s\in\Omega:\Re L_q(s)>-1,\ q=1,\ldots,n\}
\]
is nonempty. Then
\[
 M(s)=\int_{(0,1)^n}A(t;s)\prod_qt_q^{L_q(s)}\,\dd t,
 \qquad s\in\Omega_+,
\]
continues meromorphically to $\Omega$ with possible poles on $L_q+j+1=0$, $j\ge0$. For every compact $K\subset\Omega$, a finite product $H_K$ of these $L_q+j+1$, independent of $A$, makes $H_KM$ holomorphic on a neighbourhood of $K$ and bounded by finitely many $t$-derivative seminorms of $A$.
\end{lemma}

\begin{proof}
Fix $J=(J_1,\ldots,J_n)\in\mathbb Z_{>0}^n$. The multivariable Taylor formula with remainders gives
\[
 A(t;s)=\sum_{I\subseteq\{1,\ldots,n\}}\;
 \sum_{\substack{0\le\alpha_i<J_i\\ i\in I}}
 \left(\prod_{i\in I}t_i^{\alpha_i}\right)
 \left(\prod_{q\in I^c}t_q^{J_q}\right)
 B_{I,\alpha}(t_{I^c};s),
\]
where $I^c=\{1,\ldots,n\}\setminus I$. If $I^c=\varnothing$, then $B_{I,\alpha}(s)=\partial^\alpha A(0;s)/\alpha!$. $B_{I,\alpha}(s)\in C^\infty([0,1]^{I^c})$ is given by integral remainder, hence depends holomorphically on $s$ and bounded by a constant depending only on $J$ and $\|A(s)\|_{C^{J_1+\cdots+J_n}}$.

Integrating the Taylor monomials gives, initially on $\Omega_+$,
\[
 M(s)=\sum_{I,\alpha}
 \frac{\displaystyle\int_{(0,1)^{I^c}}
 B_{I,\alpha}(t_{I^c};s)
 \prod_{q\notin I}t_q^{L_q(s)+J_q}\,\dd t_{I^c}}
 {\displaystyle\prod_{i\in I}(L_i(s)+\alpha_i+1)}.
\]
For $I^c=\varnothing$, the numerator means $B_{I,\alpha}(s)$. All numerator integrals are holomorphic on the region in $\Omega$ where $\Re L_q(s)+J_q>-1$ for every $q$. Increasing $J_q$ iteratively does not change the expression on its previous domain, giving the required meromorphic continuation.
\end{proof}

\begin{proof}[Proof of Theorem~\ref{thm:Nash-continuation}]
We first continue the compactly supported distributions
\[
 \langle K_s,\varphi\rangle
 =\int_Cm(x)\varphi(x)\prod_{i=1}^r x_i^{s_i}\,\dd x,
 \qquad \varphi\in C^\infty(\R^n),\quad \Re s_i>-1.
\]
The operator in the theorem is then recovered from
$(T(s)\Phi)(b)=\langle K_s,\Phi(b,\cdot)\rangle$.
To apply boundary Taylor expansion, we first change variables so that the coefficient is smooth up to the boundary and the complex powers are monomials.

\emph{The graph and its algebraic model.}
Let $\Gamma\subset\R^{n+2}$ be the graph of $(\Re m,\Im m)$ over $C$, and define
\[
 I(\Gamma)=\{f\in\R[x_1,\ldots,x_n,u,v]:f|_\Gamma=0\},
 \qquad
 A=\R[x_1,\ldots,x_n,u,v]/I(\Gamma).
\]
We take the Zariski closure of $\Gamma$ to mean the reduced closed subscheme
$\mathbf Z=\operatorname{Spec}A\subset\mathbb A^{n+2}_{\R}$.
Embeds $A\otimes_\R\C$ into the ring of complex-valued real analytic functions on $C$. Since $C$ is connected, this ring is an integral domain. Thus $\mathbf Z$ is geometrically integral. Moreover, $\dim A=n$ by \cite[Proposition~8.4.1(i)]{BCR98}, applied to the connected Nash manifold $\Gamma$. Hence $\mathbf Z$ is an affine $\R$-variety of dimension $n$.
The coordinates $x_i,u,v$ are regular on $\mathbf Z$, and $u+iv$ restricts to $m$ on $\Gamma$.

Take a projective resolution $p:\mathbf Y\to\mathbf Z$ which is an isomorphism over the scheme-theoretic smooth locus $\mathbf Z_{\rm sm}$, using \cite[Theorem~36]{Kol07}. Put
\[
 \Gamma^\circ=\Gamma\cap\mathbf Z_{\rm sm}(\R),\qquad
 U=p(\R)^{-1}(\Gamma^\circ),\qquad
 C^\circ=\pr_x(\Gamma^\circ)\subset C,
\]
where $\pr_x:\mathbf Z\to\mathbb A^n_\R$ is the coordinate projection.
The set $\Gamma^\circ$ is open in $\mathbf Z_{\rm sm}(\R)$, and its complement in $\Gamma$ has dimension less than $n$. Consequently, $U$ is open semialgebraic, and $\pr_x\circ p$ identifies it Nash diffeomorphically with the full-measure open subset $C^\circ$ of $C$.
$\Gamma$ is bounded, so its Euclidean closure $\overline{\Gamma}$ is compact. Properness of $p$ then gives
\[
 \overline U\subset p(\R)^{-1}(\overline\Gamma),
 \qquad \overline U\ \text{compact}.
\]
As usual, closures and boundaries of real subsets are taken in the analytic topology.

\emph{The boundary.}
Write $x_i$ also for the pulled-back coordinates on $\mathbf Y$, and put
\[
 h=\prod_{i=1}^n x_i(1-x_i),\qquad
 \boldsymbol\Sigma=\mathbf Z\setminus\mathbf Z_{\rm sm}
 \quad\text{with its reduced closed subscheme structure}.
\]
The factors of $h$ record all faces of $C$. We have
\[
 \partial U\subset p(\R)^{-1}(\boldsymbol\Sigma(\R))\cup\{h=0\}.
\]
Indeed, let $y\in\overline U$ lie outside the right-hand side and put $z=p(\R)(y)$. Continuity gives $z\in\overline\Gamma$ and $\pr_x(z)\in\overline C=[0,1]^n$; since $h(y)\ne0$, in fact $\pr_x(z)\in C$. The graph is closed over $C$, so $z\in\Gamma$; as $z\notin\boldsymbol\Sigma(\R)$, we obtain $z\in\Gamma^\circ$ and hence $y\in U$. Since $U$ is open, $y\notin\partial U$.

Principalize $(h)(\mathcal I_{\boldsymbol\Sigma}\mathcal O_{\mathbf Y})$ by \cite[Theorem~35]{Kol07}, obtaining a projective morphism $q:\widetilde{\mathbf Y}\to\mathbf Y$ with smooth source and SNC reduced total-transform support $\mathbf D$. The ideal being principalized is unit along $U$, so $q$ is an isomorphism over $U$. For $\widetilde U=q(\R)^{-1}(U)$, properness and the preceding boundary inclusion give
\[
 \overline{\widetilde U}\ \text{compact},\qquad
 \widetilde U\cap\mathbf D(\R)=\varnothing,\qquad
 \partial\widetilde U\subset\mathbf D(\R).
\]
The selected real branches and the ambient maps fit into the commutative diagram
\[
\begin{tikzcd}[column sep=large,row sep=large]
 \widetilde U \arrow[r,"q|_{\widetilde U}","\sim"'] \arrow[d,hook]
 & U \arrow[r,"p|_U","\sim"'] \arrow[d,hook]
 & \Gamma^\circ \arrow[r,"\pr_x|_{\Gamma^\circ}","\sim"'] \arrow[d,hook]
 & C^\circ \arrow[d,hook] \\
 \widetilde{\mathbf Y}(\R) \arrow[r,"q(\R)"']
 & \mathbf Y(\R) \arrow[r,"p(\R)"']
 & \mathbf Z(\R) \arrow[r,"\pr_x(\R)"']
 & \R^n.
\end{tikzcd}
\]
The top arrows are Nash diffeomorphisms; the bottom arrows are induced by morphisms of $\R$-schemes. Thus integration over $C$ can be replaced by integration over $\widetilde U$.

\emph{Monomial integrals.}
Choose finitely many SNC coordinate boxes $V_j$, relatively compact in larger coordinate charts and covering $\overline{\widetilde U}$, and a smooth partition of unity $\zeta_j\in C_c^\infty(V_j)$ near $\overline{\widetilde U}$. In each chart $\mathbf D(\R)$ is a union of coordinate hyperplanes. Since $\partial\widetilde U\subset\mathbf D(\R)$, the intersection of $\widetilde U$ with each open coordinate orthant is both open and closed. Thus each orthant is either contained in $\widetilde U$ or disjoint from it. Retain the former and index them by $\alpha$. After reflection and rescaling, write $t=(t_1,\ldots,t_n)\in(0,1)^n$ for these coordinates and $\zeta_\alpha(t)$ for restriction of $\zeta_j$.

Let $\rho_\alpha$ be the coordinate expression of $\pr_x\circ p\circ q$ on such an orthant. Since $x_i$ divides $h$, the algebraic divisor of its pullback is supported on $\mathbf D$. Factoring along the local equations of this SNC divisor gives
\begin{equation}\label{eq:monomial-forms}
 \rho_{\alpha,i}(t)=u_{\alpha i}(t)\prod_qt_q^{a_{\alpha qi}},
 \qquad a_{\alpha qi}\in\mathbb Z_{\ge0},\quad 1\le i\le r.
\end{equation}
The units $u_{\alpha i}$ are smooth and nonvanishing on a neighbourhood of the closed cube, and positive there after shrinking that neighbourhood, since $\rho_{\alpha,i}>0$ on the open cube. Thus
\[
 \prod_{i=1}^r\rho_{\alpha,i}(t)^{s_i}
 =\prod_{i=1}^r u_{\alpha i}(t)^{s_i}
   \prod_qt_q^{\langle a_{\alpha q},s\rangle},
 \qquad a_{\alpha q}=(a_{\alpha q1},\ldots,a_{\alpha qr}).
\]

The map $\rho_\alpha$ and the pullback $\widetilde m_\alpha$ of $u+iv$ are smooth on a neighbourhood of the closed cube. On the open cube $\rho_\alpha$ is a diffeomorphism onto its image, $|J_\alpha|=\varepsilon_\alpha J_\alpha$, with $\varepsilon_\alpha\in\{\pm1\}, J_\alpha=\det D\rho_\alpha$.

Change of variables gives, initially for $\Re s_i>0$,
\begin{equation}\label{eq:monomial-integral}
 \langle K_s,\varphi\rangle
 =\sum_\alpha\int_{(0,1)^n}A_\alpha(t;s,\varphi)
       \prod_qt_q^{\langle a_{\alpha q},s\rangle}\,\dd t,
\end{equation}
where
\[
 A_\alpha(t;s,\varphi)=\zeta_\alpha(t)\widetilde m_\alpha(t)
     \varepsilon_\alpha J_\alpha(t)\prod_{i=1}^r u_{\alpha i}(t)^{s_i}
     \varphi(\rho_\alpha(t)),
\]
which are smooth on the closed cube, entire in $s$, continuous and linear in $\varphi\in C^\infty(\R^n)$.

\emph{Continuation and Schwartz estimates.}
Lemma~\ref{lem:Mellin} applied to \eqref{eq:monomial-integral} gives the meromorphic continuation of $K_s$, with the locally finite set of candidate poles \eqref{eq:Nash-poles} arising from the nonzero vectors $a_{\alpha q}$. For each compact set $K$ of $s$ variable, it gives a common polynomial $H_K$, independent of $\varphi$, such that $U_s=H_K(s)K_s$ is holomorphic near $K$. Also
\[
 \sup_{s\in K}|\langle U_s,\varphi\rangle|
 \le C_K\max_{|\eta|\le M}\sup_{x\in[0,1]^n}|\partial_x^\eta\varphi(x)|,
\]
with $M$ and $C_K$ independent of $\varphi$.

Define $(H_K(s)T(s)\Phi)(b)=\langle U_s,\Phi(b,\cdot)\rangle$. For affine $B$, apply the preceding estimate to $P_b\Phi(b,\cdot)$, where $P$ is any Nash differential operator on $B$. Since $P_b$ commutes with $U_s$, we obtain
\[
 \sup_{s\in K,\ b\in B}|P_b(H_K(s)T(s)\Phi)(b)|
 \le C_K\max_{|\eta|\le M}\sup_{b\in B,\ x\in[0,1]^n}
       |P_b\partial_x^\eta\Phi(b,x)|.
\]
These seminorms give \eqref{eq:Mellin-estimate}; a finite affine Nash cover and a tempered partition of unity \cite[Theorem~5.2.1]{AG08} give the estimate for general $B$. The estimates, uniform on compact sets of $s$ variable, give Cauchy bounds for the power series coefficients in every Schwartz seminorm and hence $\mcS(B)$-valued holomorphy.
\end{proof}

\section{The zero-order extension at good parameters}\label{sec:zero-order}

Throughout this section $\lambda$ is good. We glue the local regularized pairings into a distribution on $\tilde Z_w$, push it forward to $X$, and obtain $\E_{0,w,\lambda}$ (Section~\ref{ssec:gluing}). Proposition~\ref{prop:equivariance} shows that $\E_{0,w,\lambda}$ is $\fn$-equivariant, and Proposition~\ref{prop:right-inverse} that it restricts to the pushforward from $O_w$. With Proposition~\ref{prop:order-zero}, these give Theorem~\ref{thm:supported-extension} at good parameters (Theorem~\ref{thm:B-good}).

\subsection{Gluing}\label{ssec:gluing}

Fix $\gamma$ and a localized input $u_\gamma$ as in Definition~\ref{def:local-pairing}. Formula \eqref{eq:local-pairing} defines a distribution in $\mcS'(B_\gamma\times V_x,p_Z^*\mcE_\lambda)$ independent of the choice of basis of $\mcP_k(F_w)$. We first glue over the fibre.

\begin{lemma}[independence of the chart]\label{lem:overlap}
Fix $\gamma$ and two adapted charts
\[
 (\cC_v,\V_x,V_x,\boldsymbol\phi_x)
 \quad\text{and}\quad
 (\cC_{v'},\V_{x'},V_{x'},\boldsymbol\phi_{x'}).
\]
For every $\Phi\in\mcS(B_\gamma\times(V_x\cap V_{x'}),p_Z^*\mcE_\lambda)$, we have
\[
 \Bigl\langle\int^{\lambda,\gamma,v}u_\gamma,\Phi\Bigr\rangle
 =\Bigl\langle\int^{\lambda,\gamma,v'}u_\gamma,\Phi\Bigr\rangle.
\]
\end{lemma}

\begin{proof}
Write $\rho=\rho_{\sigma,\lambda}$. By \eqref{eq:eps-transition} and \eqref{eq:eps-tau},
\[
 \epsilon_{\gamma,v}=\epsilon_{\gamma,v'}\rho(R^{\gamma,v'}_{\gamma,v}),\qquad
 A_v=A_{v'}\rho(R^{\gamma,v'}_{\gamma,v}).
\]
Thus, if $\Phi=\epsilon_{\gamma,v}\boldsymbol\Phi=\epsilon_{\gamma,v'}\boldsymbol\Phi'$, then
\[
 \boldsymbol\Phi'=\rho(R^{\gamma,v'}_{\gamma,v})\boldsymbol\Phi,\qquad
 A_{v'}\boldsymbol\Phi'=A_v\boldsymbol\Phi.
\]
Hence the unregularized vector integrands agree. It remains to compare the regularizing factors. Locally on an overlap, the equations of the same components satisfy $\phi'_l=a_l\phi_l$ with nowhere vanishing regular units $a_l$. Near $z=0$ write
\[
 \prod_l|a_l|^{z_l}=1+\sum_lz_l\,h_l(z,y).
\]
After localizing a test section inside the overlap, the functions $h_l$ are smooth up to the boundary and holomorphic in $z$, with locally uniform bounds for their derivatives. Inserting them changes only the smooth amplitude in the monomial model, so the integrals below are holomorphic at $z=0$, as $\lambda$ is good. For each $p$, the difference between the two regularizations is therefore
\[
 \sum_l\Bigl\langle\mathbf u_p,\Bigl[z_l\int_{V_x\cap V_{x'}\cap F_w}A_v\boldsymbol\Phi\,h_l(z,y)\,\omega_p\prod_{l'}|\phi_{l'}|^{z_{l'}}\Bigr]_{z=0}\Bigr\rangle=0.
\]
\end{proof}

Thus the local values $\int^{\lambda,\gamma,v}u_\gamma$ glue to
$\int^{\lambda,\gamma}u_\gamma
 \in\mcS'\bigl(\tilde Z_w|_{B_\gamma},p_Z^*\mcE_\lambda\bigr)$,
supported over $K_\gamma$.

We now return to the full input $u=\sum_\gamma u_\gamma$. Using the functions $\eta_\gamma$ of Lemma~\ref{lem:coefficient-extension}, which equal one near $K_\gamma$, are supported in $B_\gamma$, and are pulled back from the base of $\tilde Z_w$, we define
\begin{equation}\label{eq:int-lambda}
 \Bigl\langle\int^\lambda u,\Psi\Bigr\rangle=\sum_\gamma\Bigl\langle\int^{\lambda,\gamma}u_\gamma,(\eta_\gamma\Psi)|_{\tilde Z_w|_{B_\gamma}}\Bigr\rangle.
\end{equation}
This defines a distribution $\int^\lambda u\in\mcS'(\tilde Z_w,p_Z^*\mcE_\lambda)$.

To obtain a distribution on $X$, we use the proper map $p_Z:\tilde Z_w\to X$, whose image is $\overline{O_w}$ by \eqref{eq:pZ-open-preimage}. Proper pullback of Schwartz functions \cite[Theorem~1.2.4(ii)]{dC91}, applied in local tempered frames, gives a continuous map
\[
 p_Z^*:\mcS(X,\mcE_\lambda)\longrightarrow\mcS(\tilde Z_w,p_Z^*\mcE_\lambda).
\]
Its strong transpose defines the proper pushforward of tempered distributions, and we put
\begin{equation}\label{eq:E0def}
\begin{aligned}
 \E_{0,w,\lambda}(u)&=(p_Z)_*\Bigl(\int^\lambda u\Bigr)
 \in\mcS'_{\overline{O_w}}(X,\mcE_\lambda),\\
 \langle\E_{0,w,\lambda}(u),\phi\rangle
 &=\Bigl\langle\int^\lambda u,p_Z^*\phi\Bigr\rangle.
\end{aligned}
\end{equation}

\subsection{Equivariance}\label{ssec:equivariance}

For $\xi\in\fn$, recall from Lemma~\ref{lem:tau-invariant} the scalar action in the product chart,
\[
 (\xi^\star f)(b,y)=\frac{\dd}{\dd t}\Big|_{t=0}f\bigl(b,\exp(-t\Ad(s_\gamma(b)^{-1})\xi)\cdot y\bigr);
\]
it differentiates only in the fibre variable. We first express the infinitesimal action in the frames $\epsilon_{\gamma,v}$.

\begin{lemma}[action in the resolution frames]\label{lem:resolution-frame-action}
Let $\Phi=\epsilon_{\gamma,v}\boldsymbol\Phi$, and write $[\Psi]_\epsilon$ for the coefficient column of a section $\Psi$ in this frame. On $B_\gamma\times C_v^\circ$,
\[
 [\xi\Phi]_\epsilon=A_v^{-1}\xi^\star(A_v\boldsymbol\Phi),\qquad
 [\xi^t\Phi]_\epsilon=-A_v^{-1}\xi^\star(A_v\boldsymbol\Phi),
\]
where $\xi^\star$ acts componentwise on columns.
\end{lemma}

\begin{proof}
The orbit frame is $N$-invariant (Lemma~\ref{lem:tau-invariant}) and $\epsilon_{\gamma,v}=\tau_\gamma A_v$ by \eqref{eq:eps-tau}. Hence
\[
 \xi\Phi=\xi(\tau_\gamma A_v\boldsymbol\Phi)
 =\tau_\gamma\xi^\star(A_v\boldsymbol\Phi)
 =\epsilon_{\gamma,v}A_v^{-1}\xi^\star(A_v\boldsymbol\Phi).
\]
\end{proof}

\begin{proposition}[equivariance]\label{prop:equivariance}
For good $\lambda$, the zero-order extension map $\E_{0,w,\lambda}$ is $\fn$-equivariant and continuous on each fixed $\IN^k$-stage. In particular, for every $k\ge0$ it gives maps
\[
 \E^k_{0,w,\lambda}:\mcS'(O_w,\mcE_\lambda)^{\IN^k}\to\mcS'_{\overline{O_w}}(X,\mcE_\lambda)^{\IN^k},
\]
which are compatible as $k$ grows.
\end{proposition}

\begin{proof}
For continuity on a fixed $\IN^k$-stage: by \eqref{eq:Mellin-estimate}, the regularized fibre operators at $z=0$ are continuous on test sections; and their strong transposes, the finite localizations and $(p_Z)_*$ are continuous. Their composite is $\E_{0,w,\lambda}$.

We differentiate the regularizing factors. It suffices to consider one localized input $u_\gamma=\tau_\gamma^\vee\sum_p(\mathbf u_p\otimes\omega_p)$ and a Schwartz test section $\Phi=\epsilon_{\gamma,v}\boldsymbol\Phi\in\mcS(B_\gamma\times V_x,p_Z^*\mcE_\lambda)$ in an adapted chart. Since $\U$ is connected, it preserves each boundary component, so its algebraic vector fields preserve the local ideals $(\phi_l)$. Hence
$\frac{\xi^\star\phi_l}{\phi_l}$
extends across $D_l$ as a regular function of the fibre coordinates, with Nash dependence on $b$.

Take $\Re z_l$ large enough to justify integration by parts, and then continue the resulting identity meromorphically. We compute
\begin{equation}\label{eq:regularized-equivariance}
\begin{aligned}
 \Bigl\langle\int^{\lambda,\gamma}(\xi u_\gamma),\Phi\Bigr\rangle
 &=-\sum_p\Bigl\langle\mathbf u_p,\Bigl[\int_{V_x\cap F_w}\xi^\star\Bigl(A_v\boldsymbol\Phi\prod_l|\phi_l|^{z_l}\Bigr)\omega_p\Bigr]_{z=0}\Bigr\rangle\\
 &=-\sum_p\Bigl\langle\mathbf u_p,\Bigl[\int_{V_x\cap F_w}\xi^\star(A_v\boldsymbol\Phi)\,\omega_p\prod_l|\phi_l|^{z_l}\Bigr]_{z=0}\Bigr\rangle\\
 &\quad-\sum_{p,l}\Bigl\langle\mathbf u_p,\Bigl[z_l\int_{V_x\cap F_w}A_v\boldsymbol\Phi\,\frac{\xi^\star\phi_l}{\phi_l}\,\omega_p\prod_{l'}|\phi_{l'}|^{z_{l'}}\Bigr]_{z=0}\Bigr\rangle.
\end{aligned}
\end{equation}

The last term is zero since $\lambda$ is good, the integral is holomorphic at $z=0$, and it is multiplied by $z_l$. Lemma~\ref{lem:resolution-frame-action} gives
\[
 A_v[\xi^t\Phi]_\epsilon=-\xi^\star(A_v\boldsymbol\Phi),
\]
and hence
\[
 \Bigl\langle\int^{\lambda,\gamma}(\xi u_\gamma),\Phi\Bigr\rangle=\Bigl\langle\int^{\lambda,\gamma}u_\gamma,\xi^t\Phi\Bigr\rangle.
\]
The chartwise equalities give infinitesimal equivariance on $\tilde Z_w|_{B_\gamma}$. Both the input partition $\beta_\gamma$ and the extension cutoff $\eta_\gamma$ depend only on the base, so they commute with $\xi$. Summing over $\gamma$ and using the $N$-equivariance of $p_Z$, we obtain
\[
 \E_{0,w,\lambda}(\xi u)=\xi\,\E_{0,w,\lambda}(u).
\]
\end{proof}

\subsection{Restriction to the orbit}\label{ssec:right-inverse}

\begin{proposition}\label{prop:right-inverse}
Let $i:O_w\hookrightarrow\mcU_w$ be the closed embedding and $r^X_w$ the restriction from $X$ to $\mcU_w$. Then
\[
 r^X_w\circ\E^k_{0,w,\lambda}=i_*\qquad\text{on }\ \mcS'(O_w,\mcE_\lambda)^{\IN^k}.
\]
\end{proposition}

\begin{proof}
Let $\phi\in\mcS(\mcU_w,\mcE_\lambda)$ and let $\widetilde\phi$ be its extension by zero to $X$. By rapid decay of $i^*\phi$, evaluation of the local regularized integrals at $z=0$ is ordinary integration. We obtain
\[
\begin{aligned}
 \langle r_w^X\E_{0,w,\lambda}(u),\phi\rangle
 =\Bigl\langle\int^\lambda u,p_Z^*\widetilde\phi\Bigr\rangle
 =\sum_\gamma\langle u_\gamma,i^*\phi\rangle
 =\langle i_*u,\phi\rangle.
\end{aligned}
\]
\end{proof}

\subsection{Theorem~\ref{thm:supported-extension} at good parameters}\label{ssec:thmB-good}

\begin{theorem}\label{thm:B-good}
Let $\lambda$ be good. Then $\E_0=\E_{0,w,\lambda}$ satisfies the hypotheses of Proposition~\ref{prop:order-zero}, and the resulting map $\E_{w,\lambda}:=\E_{w,V_{\sigma,\lambda}}$ is a linear section of the restriction map
\[
 \mcS'_{\overline{O_w}}(X,\mcE_\lambda)^{[\fn]}\to\mcS'_{O_w}(\mcU_w,\mcE_\lambda)^{[\fn]}.
\]
In particular, Theorem~\ref{thm:supported-extension} holds for $O_w$ at $\lambda$, for every Nash open subset $U_w$ in which $O_w$ is closed.
\end{theorem}

\begin{proof}
By Propositions~\ref{prop:equivariance} and~\ref{prop:right-inverse}, $\E_{0,w,\lambda}$ maps $\mcS'(O_w,\mcE_\lambda)^{[\fn]}$ to $\mcS'_{\overline{O_w}}(X,\mcE_\lambda)^{[\fn]}$ and satisfies $r^X_w\circ\E_{0,w,\lambda}=i_*$. Proposition~\ref{prop:order-zero} gives the linear section. Since the image of $\E_{w,\lambda}$ is supported in $\overline{O_w}$, Lemma~\ref{lem:supported-nbhd} transfers the conclusion to any $U_w$ in which $O_w$ is closed.
\end{proof}

\section{Bad parameters}\label{sec:bad}

We complete the proof of Theorem~\ref{thm:supported-extension} by the deformation and pole-removal argument of \cite[Lemmas~4.5--4.6]{Abe12}.

\subsection{Meromorphic deformation}\label{ssec:deformation}

Fix $\lambda_0\in\fa_{0,\C}^*$. By Remark~\ref{rem:transverse}, the line $\lambda(\tau)=\lambda_0+\tau\nu_0$ consists of good parameters for every $Q$-orbit when $0<|\tau|<\epsilon$, for some $\epsilon>0$. Use the compact picture
\[
 C_\lambda:\mcE_K=K\times^{M_0}V_\sigma\xrightarrow{\sim}\mcE_\lambda,
 \qquad [k,v]\longmapsto[k,v].
\]
For any Nash open $U\subset X$, holomorphy and Laurent expansions below are understood after transport by $C_{\lambda(\tau)}^*$ to the fixed strong dual $\mcS'(U,\mcE_K)$; finite pole order means that multiplication by some $\tau^p$ gives a holomorphic family. A family is \emph{uniformly $\IN$-finite} if one power of $\IN$ annihilates it for all parameters under consideration.

We recall the parameter dependence in this picture. For the dual orbit frame of \eqref{eq:tau-frame}, with $g=g_\gamma(b,f)$ and $H(g)=\log a(g)$,
\begin{equation}\label{eq:tau-compact}
 C_\lambda^*\tau_{\gamma,\lambda}^\vee
 =e^{-\lambda(H(g))}[\kappa(g),\mathbf e^\vee]_K.
\end{equation}
The frame changes in \eqref{eq:tau-compact} act holomorphically on Schwartz spaces and their strong duals. The transported action $A_D(\tau)$ of $D\in\mcU(\fg_\C)$ also depends holomorphically on $\tau$, by the compact-picture calculation in \cite[p.~444]{Abe12}.

\begin{proposition}[meromorphic deformation]\label{prop:meromorphic-deformation}
For every $u_0\in\mcS'_{O_w}(\mcU_w,\mcE_{\lambda_0})^{[\fn]}$ there is a holomorphic family $u(\tau)\in\mcS'_{O_w}(\mcU_w,\mcE_{\lambda(\tau)})$ with $u(0)=u_0$ and a meromorphic family $S(\tau)\in\mcS'(X,\mcE_{\lambda(\tau)})$ of finite pole order at $0$ such that
\[
 S(\tau)|_{\mcU_w}=u(\tau),\qquad
 \supp S(\tau)\subset\overline{O_w}.
\]
Both families are uniformly $\IN$-finite.
\end{proposition}

\begin{proof}
By Lemma~\ref{lem:Borel}, $u_0$ has finite normal order $r$. Remark~\ref{rem:section-sr} and Proposition~\ref{prop:poly-densities} give a finite presentation
\[
 u_0=\sum_{\gamma,i}D_{\gamma i}\widetilde{v_{\gamma i}(0)},\qquad
 v_{\gamma i}(0)=\tau_{\gamma,\lambda_0}^\vee\mathbf v_{\gamma i},\qquad
 \mathbf v_{\gamma i}\in
 \bigl(\mcS'(B_\gamma)\otimes\mcP_k(F_w)\bigr)^d,
\]
where $D_{\gamma i}\in\mcU_{\le r}(\fg_\C)$, the coefficient columns are supported in $C_\gamma$, and the tilde denotes the fixed support extension of Lemma~\ref{lem:coefficient-extension}. Keep these columns fixed and put
\[
 v_{\gamma i}(\tau)=\tau_{\gamma,\lambda(\tau)}^\vee\mathbf v_{\gamma i},\qquad
 u(\tau)=\sum_{\gamma,i}D_{\gamma i}\widetilde{v_{\gamma i}(\tau)}.
\]
For $0<|\tau|<\epsilon$, define
\[
 S(\tau)=\sum_{\gamma,i}D_{\gamma i}
 \E_{0,w,\lambda(\tau)}\bigl(\widetilde{v_{\gamma i}(\tau)}\bigr).
\]
The family $u(\tau)$ is holomorphic by \eqref{eq:tau-compact}. The orbit frames are $N$-invariant (Lemma~\ref{lem:tau-invariant}), and Proposition~\ref{prop:equivariance} preserves the annihilating power $\IN^k$ of the coefficients. Lemma~\ref{lem:commutator} therefore annihilates both families by $\IN^{k+b_r-1}$. Proposition~\ref{prop:right-inverse} gives the restriction identity, and the extension construction gives the support assertion.

It remains to prove that $S(\tau)$ has finite pole order. Localize each $\widetilde{v_{\gamma i}(\tau)}$ using the fixed cutoffs in Definition~\ref{def:local-pairing}; write $\gamma'$ for the chart index in this second localization. These localizations and the changes of orbit frame preserve $\IN^k$-annihilation. By Proposition~\ref{prop:poly-densities}, the resulting coefficient columns lie in $\mcS'(B_{\gamma'})\otimes\mcP_k(F_w)\otimes V_\sigma^*$. In a fixed basis of $\mcP_k(F_w)$, \eqref{eq:local-pairing} gives a finite sum
\[
 \bigl\langle C_{\lambda(\tau)}^*S(\tau),\Phi\bigr\rangle
 =\sum_j\bigl\langle\mathbf b_j(\tau),
 J_{\omega_j}\bigl(0,\lambda(\tau);L_j(\tau)\Phi\bigr)\bigr\rangle,
 \qquad \Phi\in\mcS(X,\mcE_K).
\]
Here $j$ collects the chart, coefficient and differential-operator indices, with base-chart index $\gamma'_j$; $\mathbf b_j(\tau)\in\mcS'(B_{\gamma'_j})^d$ is holomorphic and $\omega_j\in\mcP_k(F_w)$ is fixed. For a term coming from $D_{\gamma i}$, the test-section map has the form
\[
 L_j(\tau)\Phi
 =\Bigl[\chi_j\,p_Z^*\bigl(D_{\gamma i}^{t}
       (C_{\lambda(\tau)}\Phi)\bigr)\Bigr]_{\epsilon_{\gamma'_j,v_j}}.
\]
Here $\chi_j$ collects the fixed cutoffs, and the brackets denote the coefficient column in the resolution frame at $\lambda(\tau)$ on $B_{\gamma'_j}\times V_{x_j}$. Thus $L_j(\tau):\mcS(X,\mcE_K)\to\mcS(B_{\gamma'_j}\times V_{x_j})^d$ is a holomorphic family of continuous maps.

For each term, Theorem~\ref{thm:Nash-continuation} supplies a finite product $H_j(z,\tau)$ of factors
\[
 \langle a_\nu,z\rangle+\ell_\nu(\lambda(\tau))+N,
 \qquad N\in\mathbb Z,
\]
which clears its poles near $(z,\tau)=(0,0)$, independently of $\Phi$; here $z$ is the regularizing variable for the corresponding chart. At $z=0$, each factor is $\ell_\nu(\lambda_0)+N+\tau\ell_\nu(\nu_0)$. By Remark~\ref{rem:transverse}, its coefficient of $\tau$ is nonzero. Consequently, after shrinking the parameter neighbourhood,
\[
 H(\tau):=\prod_jH_j(0,\tau)=\tau^p h(\tau),
 \qquad h(0)\ne0.
\]
This gives holomorphy of $\tau^pC_{\lambda(\tau)}^*S(\tau)$.
\end{proof}

\subsection{Pole removal and completion of the proof}\label{ssec:pole-removal}

\begin{proposition}[holomorphic extension]\label{prop:holomorphic-extension}
Put $X_w^\circ=X\setminus(\overline{O_w}\setminus O_w)$. For every $u_0\in\mcS'_{O_w}(X_w^\circ,\mcE_{\lambda_0})^{[\fn]}$, there is, near $\tau=0$, a holomorphic, uniformly $\IN$-finite family
\[
 F(\tau)\in\mcS'_{\overline{O_w}}(X,\mcE_{\lambda(\tau)}),
 \qquad F(0)|_{X_w^\circ}=u_0.
\]
\end{proposition}

\begin{proof}
We follow the induction in the proof of \cite[Lemma~4.6, p.~445]{Abe12}. Choose a filtration
\[
 \varnothing=Z_0\subset Z_1\subset\cdots\subset Z_m=\overline{O_w}
\]
by closed unions of $Q$-orbits, with $O_j=Z_j\setminus Z_{j-1}$ a single orbit and $O_m=O_w$. Along this filtration we also prove that every $v_0\in\mcS'_{Z_j}(X,\mcE_{\lambda_0})^{[\fn]}$ admits a holomorphic, uniformly $\IN$-finite family $v(\tau)$ supported in $Z_j$, with $v(0)=v_0$. This assertion is trivial for $Z_0$.

Assume the assertion for $Z_{j-1}$ and put $V_j=X\setminus Z_{j-1}$. For $u_0\in\mcS'_{O_j}(V_j,\mcE_{\lambda_0})^{[\fn]}$, Proposition~\ref{prop:meromorphic-deformation} and Lemma~\ref{lem:supported-nbhd} give a uniformly $\IN$-finite meromorphic family $S(\tau)$ supported in $Z_j$, whose restriction to $V_j$ is holomorphic and takes the value $u_0$ at $0$. If it has a pole, write in the fixed compact picture
\[
 C_{\lambda(\tau)}^*S(\tau)=\sum_{q=-p}^{\infty}\tau^qS_q,
 \qquad p>0.
\]
Holomorphy on $V_j$ gives $\supp S_{-p}\subset Z_{j-1}$. If $\IN^kS(\tau)=0$, the lowest Laurent coefficient of $A_D(\tau)C_{\lambda(\tau)}^*S(\tau)=0$ gives
\[
 A_D(0)S_{-p}=0\qquad(D\in\IN^k).
\]
By the inductive hypothesis, there is a holomorphic, uniformly $\IN$-finite family $L(\tau)$ supported in $Z_{j-1}$, with $C_{\lambda_0}^*L(0)=S_{-p}$. Thus
\[
 S(\tau)-\tau^{-p}L(\tau)
\]
has smaller pole order and the same restriction to $V_j$. Repeating finitely many times gives a holomorphic family $F(\tau)$ supported in $Z_j$, with $F(0)|_{V_j}=u_0$. At each subtraction, a common annihilating power of $\IN$ is obtained by taking the larger of the two exponents, so $F(\tau)$ remains uniformly $\IN$-finite.

To complete the induction, take $v_0\in\mcS'_{Z_j}(X,\mcE_{\lambda_0})^{[\fn]}$ and apply this construction to $v_0|_{V_j}$. The difference $v_0-F(0)$ is supported in $Z_{j-1}$. By induction it has a holomorphic, uniformly $\IN$-finite deformation supported there; adding this deformation to $F(\tau)$ proves the assertion for $Z_j$.

Taking $j=m$ gives the required family, since $V_m=X_w^\circ$.
\end{proof}

\begin{proof}[Proof of Theorem~\ref{thm:supported-extension}]
Evaluation at $\tau=0$ in Proposition~\ref{prop:holomorphic-extension} gives an $\IN$-finite extension supported in $\overline{O_w}$. Lemma~\ref{lem:supported-nbhd} transfers this conclusion to any neighbourhood in which $O_w$ is closed, proving Theorem~\ref{thm:supported-extension} at $\lambda_0$.
\end{proof}

\section{General parabolic subgroups}\label{sec:general-parabolics}

We now prove Theorem~\ref{thm:general-extension}. We may assume that $P=L_PU_P$ and $Q=LN$ are standard, and put $X=G/P$. For $E\in\mathrm{CW}(G)$, write
\[
 \mathrm{CJ}^{\infty}_N(E)=\varprojlim_k E/\overline{I_N^kE},
\]
where the closures are taken in the Fr\'echet topology and the limit carries the projective limit topology. We use the same notation for Levi subgroups and their parabolic unipotent radicals.

For $F\in\mathrm{CW}(L_P)$, denote its action by
$\sigma_F$. Put
\[
 J_W(F)=\mcS'(W,\mcE_F)^{[\fn]},\qquad \mcE_F=G\times^P F
\]
for every $Q$-stable Nash open subset $W\subset X$.

\subsection{Reduction to coefficient lifting}\label{ssec:coefficient-reduction}

We reduce Theorem~\ref{thm:general-extension}(i) to coefficient lifting with support on an orbit. First, induction on the length of $F$ reduces the assertion to irreducible coefficients. For reducible $F$, choose an irreducible subrepresentation $F_1\subset F$. Exactness of parabolic induction and strict exactness of the smooth Casselman--Jacquet functor \cite[Theorem~14.3]{CWYZ26} give, by continuous duality, the upper exact row in
\[
\begin{tikzcd}[column sep=small,row sep=large]
 0\arrow[r]&J_X(F/F_1)\arrow[r]\arrow[d,two heads,"\mathrm{res}"']
 &J_X(F)\arrow[r]\arrow[d,"\mathrm{res}"]
 &J_X(F_1)\arrow[r]\arrow[d,two heads,"\mathrm{res}"]&0\\
 0\arrow[r]&J_W(F/F_1)\arrow[r]&J_W(F)\arrow[r]&J_W(F_1).
\end{tikzcd}
\]
The lower row is left exact, and the outer vertical maps are surjective by induction, so the middle map is surjective.

Now suppose that $F$ is irreducible. Put $P_{0,P}=P_0\cap L_P$, $X_0=G/P_0$, and let $p:X_0\to X$ be the projection. Casselman's subrepresentation theorem gives a closed embedding
\[
 j:F\hookrightarrow\Pi_\rho:=\operatorname{Ind}_{P_{0,P}}^{L_P}\rho,
 \qquad \rho=\eta\otimes e^\mu,
\]
with $\eta$ a finite-dimensional representation of $M_0$ \cite[Corollary~5.7]{Cas89}. Write $j_W':J_W(\Pi_\rho)\to J_W(F)$ for the induced map, and put $\mathcal F_\rho=G\times^{P_0}\rho$. All inductions are unnormalized. Induction in stages gives the topological isomorphism
\begin{equation}\label{eq:open-induction-stages}
 \mcS(W,\mcE_{\Pi_\rho})\simeq\mcS(p^{-1}W,\mathcal F_\rho).
\end{equation}
It commutes with restriction and the $\fn$-action, giving the diagram
\begin{equation}\label{eq:principal-coefficient-extension}
\begin{tikzcd}[column sep=large,row sep=large]
 \mcS'(X_0,\mathcal F_\rho)^{[\fn]}\arrow[r,"\sim"]\arrow[d,two heads,"\mathrm{res}"']
 &J_X(\Pi_\rho)\arrow[r,"j_X'"]\arrow[d,two heads,"\mathrm{res}"]
 &J_X(F)\arrow[d,"\mathrm{res}"]\\
 \mcS'(p^{-1}W,\mathcal F_\rho)^{[\fn]}\arrow[r,"\sim"]
 &J_W(\Pi_\rho)\arrow[r,"j_W'"']&J_W(F).
\end{tikzcd}
\end{equation}
The left vertical map is surjective by Theorem~\ref{thm:minimal-extension}, so $\Pi_\rho$ has the extension property for every $Q$-stable Nash open subset of $X$. It therefore suffices to prove that $j_W'$ is surjective.

Induct on the number of $Q$-orbits in $W$, the empty case being immediate. Choose an orbit $O$ closed in $W$, put $W'=W\setminus O$, and consider
\[
\begin{tikzcd}[column sep=large,row sep=large]
 J_X(\Pi_\rho)\arrow[r,"j_X'"]\arrow[d,two heads,"\mathrm{res}"']
 &J_X(F)\arrow[d,"\mathrm{res}"]\\
 J_W(\Pi_\rho)\arrow[r,"j_W'"]\arrow[d,"\mathrm{res}"']
 &J_W(F)\arrow[d,"\mathrm{res}"]\\
 J_{W'}(\Pi_\rho)\arrow[r,two heads,"j_{W'}'"']&J_{W'}(F).
\end{tikzcd}
\]
For $T\in J_W(F)$, lift $T|_{W'}$ through the bottom map by induction and extend the lift to $S\in J_X(\Pi_\rho)$. The difference $T-j_W'(S|_W)$ is supported in $O$. Lifting this difference and adding the lift to $S|_W$ completes the induction. Thus the remaining assertion is the following, and the proof is given in Section~\ref{ssec:orbit-lifting}:

\begin{proposition}[coefficient lifting with support]\label{prop:coefficient-supported}
Let $F_1\hookrightarrow F_2$ be a Casselman--Wallach embedding, and let $O$ be a $Q$-orbit closed in a Nash open subset $U\subset X$. Then
\[
 \mcS'_O(U,\mcE_{F_2})^{[\fn]}\longrightarrow
 \mcS'_O(U,\mcE_{F_1})^{[\fn]}
\]
is surjective. A distribution of normal order at most $r$ has a lift of normal order at most $r$.
\end{proposition}

Assuming Proposition~\ref{prop:coefficient-supported}, the preceding argument proves Theorem~\ref{thm:general-extension}(i). For part~(ii), let $F_1\hookrightarrow F_2$ be a Casselman--Wallach embedding and $T\in J_W(F_1)$. By part~(i), extend $T$ to $\widetilde T\in J_X(F_1)$. Exactness of parabolic induction and strict exactness of the smooth Casselman--Jacquet functor give a lift $S\in J_X(F_2)$, whose restriction to $W$ lifts $T$. Together with left exactness, this proves part~(ii).

\subsection{Lifting along an orbit}\label{ssec:orbit-lifting}

Fix a $Q$-orbit $O_w\subset X$. The subgroup
\[
 P_w=L_P\cap\dot w^{-1}Q\dot w
\]
is parabolic in $L_P$, with unipotent radical $V_w=L_P\cap\dot w^{-1}N\dot w$. Put $\mathfrak v_w=\Lie V_w$ and $I_{V_w}=\ker(\mcU((\mathfrak v_w)_\C)\to\C)$.

Next lemma is immediate by the strict exactness of $\mathrm{CJ}^{\infty}_{V_w}$ and nuclearity of $\mcS(B)$.
\begin{lemma}\label{lem:passive-jacquet}
Let $F_1\hookrightarrow F_2$ be an embedding in $\mathrm{CW}(L_P)$, and let $B$ be a Nash manifold on which $V_w$ acts trivially. Then restriction induces a surjection
\[
 \bigl(\mcS(B)\widehat\otimes F_2\bigr)'^{[\mathfrak v_w]}
 \longrightarrow
 \bigl(\mcS(B)\widehat\otimes F_1\bigr)'^{[\mathfrak v_w]}.
\]
If $t$ is supported in a compact set $K\subset B$, its lift may be chosen with support in any prescribed compact neighbourhood $K'$ of $K$ in $B$.
\end{lemma}

Use the charts of Section~\ref{ssec:orbit-geometry}. Recall that
\[
 N_w=N\cap\dot wP\dot w^{-1},\qquad
 V_w=L_P\cap\dot w^{-1}N\dot w.
\]
As in Section~\ref{ssec:orbit-frames}, write points of $\Omega_\gamma$ as $(b,y)\in B_\gamma\times F_w$, corresponding to $s_\gamma(b)y$, and use the orbit lift
\[
 g_\gamma:\B_\gamma\times\F_w\longrightarrow\G,
 \qquad g_\gamma(b,y)=s_\gamma(b)c_w(y)\dot w.
\]
The map
\[
 \begin{aligned}
 \tau_{\gamma,F}:\Omega_\gamma\times F
 &\xrightarrow{\sim}\mcE_F|_{\Omega_\gamma},\\
 ((b,y),v)&\longmapsto[g_\gamma(b,y),v]
 \end{aligned}
\]
is a tempered orbit trivialization \cite[Section~3.3]{CS21}. Writing
$\phi(b,y)=\tau_{\gamma,F}((b,y),\boldsymbol\phi_\gamma(b,y))$ gives the topological identification
\[
 \mcS(\Omega_\gamma,\mcE_F)
 \simeq\mcS(B_\gamma\times F_w)\widehat\otimes F,
 \qquad\phi\longmapsto\boldsymbol\phi_\gamma.
\]

Let $\operatorname{pr}_{L_P}:P\to L_P$ be the Levi projection, and write
\[
 \theta:N_w\twoheadrightarrow V_w,\qquad
 h\longmapsto\operatorname{pr}_{L_P}(\dot w^{-1}h\dot w).
\]
For $n\in N$, put $a_b(n)=s_\gamma(b)^{-1}ns_\gamma(b)$ and factor
\[
 a_b(n)^{-1}c_w(y)=c_w(y_n(b,y))h_n(b,y),
 \qquad y_n(b,y)\in F_w,\quad h_n(b,y)\in N_w.
\]
This factorization is regular by Remark~\ref{rem:cw}, and $y_n(b,y)=a_b(n)^{-1}y$. Put $c_n(b,y)=\theta(h_n(b,y))\in V_w$. The natural action on the coefficient function $\boldsymbol\phi_\gamma$ is
\begin{equation}\label{eq:PQ-orbit-action}
 (n\cdot\boldsymbol\phi_\gamma)(b,y)
 =\sigma_F(c_n(b,y))^{-1}\boldsymbol\phi_\gamma(b,y_n(b,y)).
\end{equation}
Writing $\xi^\star$ for its differential, we obtain, for $\xi\in\fn$,
\[
 \begin{aligned}
 (\xi^\star\boldsymbol\phi_\gamma)(b,y)
 &=\frac{\mathrm d}{\mathrm dt}\Big|_{t=0}
   \left[\sigma_F(c_{\exp(t\xi)}(b,y))^{-1}
   \boldsymbol\phi_\gamma(b,y_{\exp(t\xi)}(b,y))\right]\\
 &=(Y_\xi\boldsymbol\phi_\gamma)(b,y)
  -\mathrm d\sigma_F(\dot c_\xi(b,y))\boldsymbol\phi_\gamma(b,y),
 \end{aligned}
\]
where
\[
 \begin{aligned}
 Y_\xi(b,y)&=\frac{\mathrm d}{\mathrm dt}\Big|_{t=0}
   y_{\exp(t\xi)}(b,y)\in T_yF_w,\\
 \dot c_\xi(b,y)&=\frac{\mathrm d}{\mathrm dt}\Big|_{t=0}
   c_{\exp(t\xi)}(b,y)\in\mathfrak v_w.
 \end{aligned}
\]
Thus the action differentiates only in the fibre variable, with the additional coefficient action through $\sigma_F$.

As in Section~\ref{ssec:poly-densities}, let $\Omega$ be a nonzero $\U$-invariant algebraic top form on $\F_w$, put $\mu=|\Omega|$, and write $\mcP(F_w)=\mathcal O(\F_w)_\C\,\mu$.

\begin{proposition}[regularity along the fibre]\label{prop:PQ-poly-densities}
In the orbit trivialization $\tau_{\gamma,F}$, there is a vector-space identification, natural in $F$,
\begin{equation}\label{eq:PQ-polynomial-coefficients}
 \mcS'(\Omega_\gamma,\mcE_F)^{[\fn]} = \bigl(\mcS(B_\gamma)\widehat\otimes F\bigr)'^{[\mathfrak v_w]}
 \otimes_\C\mcP(F_w).
\end{equation}
Here $\otimes_\C$ is the algebraic tensor product, and the $V_w$-action on the dual is induced from $F$, with trivial action on $B_\gamma$. This identification preserves support over the base.
\end{proposition}
\begin{proof}
Set $H=N\cap\dot wU_P^-\dot w^{-1}$, identified with $F_w$ by $c_w$ from Remark~\ref{rem:cw}. For $X\in\fn$, let $\mathscr D_X=(\Ad(s_\gamma(b))X)^\star$. Since these operators do not differentiate $b$, $\IN^ku=0$ if and only if all products of $k$ operators $\mathscr D_X$ annihilate $u$.

For $X\in\Lie H$, substituting $n=s_\gamma(b)\exp(tX)s_\gamma(b)^{-1}$ in \eqref{eq:PQ-orbit-action} gives $y_n=\exp(-tX)y$ and $h_n=1$, since $\exp(-tX)c_w(y)\in H$. Hence
\[
 \mathscr D_X\boldsymbol\phi_\gamma(b,y)
 =\left.\frac{\mathrm d}{\mathrm dt}\right|_{t=0}
   \boldsymbol\phi_\gamma(b,\exp(-tX)y)
\]
is the natural translation action. Put $E=(\mcS(B_\gamma)\widehat\otimes F)'$. By \cite[Proposition~A.4]{Abe12}, $\IN^ku=0$ therefore gives
\[
 u\in E\otimes_\C\mcP(F_w).
\]

Write $k_v\star$ for the action of $\Ad(s_\gamma(b))k_v\cdot$ in \eqref{eq:PQ-orbit-action}. The subgroup $\dot wV_w\dot w^{-1}\subset N_w$ normalizes $H$ and maps isomorphically onto $V_w$ under $\theta$. For $k_v=\dot wv\dot w^{-1}$, the factorization $k_v^{-1}h=(k_v^{-1}hk_v)k_v^{-1}$ gives
\[
 k_v\star(\ell\otimes p\mu)
 =\sigma_E(v)\ell\otimes p(k_v^{-1}(\,\cdot\,)k_v)\mu
 \qquad(\ell\in E,\ p\in\mathcal O(\F_w)_\C).
\]
Choose a finite-dimensional rational $V_w$-submodule $M\subset\mcP(F_w)$ with $u\in E\otimes M$, and $a\ge1$ with $I_{V_w}^aM^*=0$. Since $I_{V_w}^ku=0$ for the diagonal action, Lemma~\ref{lem:commutator}(a) gives that every coefficient lies in $E^{[\mathfrak v_w]}$. The converse inclusion is immediate.
\end{proof}

Combining Proposition~\ref{prop:PQ-poly-densities} and Lemma~\ref{lem:passive-jacquet} gives the following:
\begin{lemma}\label{lem:PQ-local-lifting}
Let $F_1\hookrightarrow F_2$ be a Casselman--Wallach embedding of $L_P$. Every $\IN$-finite distribution on $\Omega_\gamma$ with coefficient $\mcE_{F_1}$, supported over $K\Subset B_\gamma$, has an $\IN$-finite lift with coefficient $\mcE_{F_2}$, supported over a compact neighbourhood of $K$ in $B_\gamma$.
\end{lemma}

\begin{proof}[Proof of Proposition~\ref{prop:coefficient-supported}]
By Lemma~\ref{lem:supported-nbhd}, we may take $U=\mcU_w$.

For $T\in F_r\mcS'_{O_w}(\mcU_w,\mcE_{F_1})^{[\fn]}$, Theorem~\ref{thm:global-presentation} and Remark~\ref{rem:section-sr} give
\[
 T=\sum_{\gamma,i}D_{\gamma i}\widetilde u_{\gamma i},\qquad
 \deg D_{\gamma i}\leq r,
\]
where $u_{\gamma i}\in\mcS'(\Omega_\gamma,\mcE_{F_1})^{[\fn]}$ has compact base support and $\widetilde u_{\gamma i}$ denotes extension to $O_w$. Lift each $u_{\gamma i}$ by Lemma~\ref{lem:PQ-local-lifting} and extend the lift to $O_w$ by Lemma~\ref{lem:coefficient-extension}. Applying the same $D_{\gamma i}$ and summing gives a lift of $T$, of normal order at most $r$. It is $\IN$-finite by Lemma~\ref{lem:commutator}.
\end{proof}


\begin{thebibliography}{99}
\bibitem{Abe12}
N.~Abe,
\emph{Generalized Jacquet modules of parabolically induced representations},
Publ. Res. Inst. Math. Sci. \textbf{48} (2012), no.~2, 419--473.

\bibitem{AG08}
A.~Aizenbud and D.~Gourevitch,
\emph{Schwartz functions on Nash manifolds},
Int. Math. Res. Not. IMRN 2008, no.~5, Art. ID rnm155, 37~pp.

\bibitem{Ati70}
M.~F.~Atiyah,
\emph{Resolution of singularities and division of distributions},
Comm. Pure Appl. Math. \textbf{23} (1970), no.~2, 145--150.

\bibitem{BM88}
E.~Bierstone and P.~D.~Milman,
\emph{Semianalytic and subanalytic sets},
Publ. Math. Inst. Hautes \'Etudes Sci. \textbf{67} (1988), 5--42.

\bibitem{BCR98}
J.~Bochnak, M.~Coste, and M.-F.~Roy,
\emph{Real Algebraic Geometry},
Ergeb. Math. Grenzgeb. (3), vol.~36,
Springer-Verlag, Berlin, 1998.

\bibitem{BT72}
A.~Borel and J.~Tits,
\emph{Compl\'ements \`a l'article ``Groupes r\'eductifs''},
Publ. Math. Inst. Hautes \'Etudes Sci. \textbf{41} (1972), 253--276.

\bibitem{Cas89}
W.~Casselman,
\emph{Canonical extensions of Harish--Chandra modules to representations of $G$},
Canad. J. Math. \textbf{41} (1989), no.~3, 385--438.

\bibitem{CHM00}
W.~Casselman, H.~Hecht, and D.~Mili\v{c}i\'{c},
\emph{Bruhat filtrations and Whittaker vectors for real groups},
in \emph{The Mathematical Legacy of Harish-Chandra} (Baltimore, MD, 1998),
Proc. Sympos. Pure Math., vol.~68,
American Mathematical Society, Providence, RI, 2000, pp.~151--190.

\bibitem{CS21}
Y.~Chen and B.~Sun,
\emph{Schwartz homologies of representations of almost linear Nash groups},
J. Funct. Anal. \textbf{280} (2021), no.~7, Art. ID 108817, 50~pp.

\bibitem{CWYZ26}
K.~Y.~Chan, K.~Wu, J.~Yu, and H.~Zhang,
\emph{Stability of the smooth Casselman--Jacquet functor},
\href{https://arxiv.org/abs/2606.06243v2}{arXiv:2606.06243v2}, 2026.

\bibitem{dC91}
F.~du~Cloux,
\emph{Sur les repr\'esentations diff\'erentiables des groupes de Lie alg\'ebriques},
Ann. Sci. \'Ec. Norm. Sup\'er. (4) \textbf{24} (1991), no.~3, 257--318.

\bibitem{Kol07}
J.~Koll\'ar,
\emph{Resolution of singularities---Seattle lecture},
\href{https://arxiv.org/abs/math/0508332v3}{arXiv:math/0508332v3}, 2007.

\bibitem{Par94}
A.~Parusi\'nski,
\emph{Subanalytic functions},
Trans. Amer. Math. Soc. \textbf{344} (1994), no.~2, 583--595.

\bibitem{Sun15}
B.~Sun,
\emph{Almost linear Nash groups},
Chinese Ann. Math. Ser. B \textbf{36} (2015), no.~3, 355--400.

\bibitem{SZ12}
B.~Sun and C.-B.~Zhu,
\emph{Multiplicity one theorems: the Archimedean case},
Ann. of Math. (2) \textbf{175} (2012), no.~1, 23--44.

\bibitem{XY26}
J.~Xiao and J.~Yu,
\emph{A proof of Casselman's comparison conjecture},
\href{https://arxiv.org/abs/2611.25367}{arXiv:2611.25367}, 2026.

\end{thebibliography}
\end{document}